\documentclass[nonblindrev]{informs3_noinforms}              

\usepackage[numbers]{natbib}

\usepackage{enumerate}
\newtheorem{observation}{Observation}

\newcommand{\pr}{\mathbb{P}}
\newcommand{\E}{\mathbb{E}}

\newcommand{\reals}{{\mathbb R}}
\newcommand{\ignore}[1]{\relax}

\newcommand{\sign}{\text{sign}}

\newcommand{\N}{{\cal N}}

\renewcommand{\theequation}{\thesection.\arabic{equation}}

\usepackage{natbib}
 \NatBibNumeric
 \def\bibfont{\small}%
 \def\bibsep{\smallskipamount}%
 \bibpunct[, ]{[}{]}{,}{n}{}{,}%

\usepackage[colorlinks=true,breaklinks=true,bookmarks=true,urlcolor=blue,
     citecolor=blue,linkcolor=blue,bookmarksopen=false,draft=false]{hyperref}

\def\EMAIL#1{\href{mailto:#1}{#1}}
\def\URL#1{\href{#1}{#1}}         

\TheoremsNumberedThrough     

\EquationsNumberedThrough    

\begin{document}



\RUNTITLE{Heavy-tailed queues in the Halfin-Whitt regime}

\TITLE{Heavy-tailed queues in the Halfin-Whitt regime}

\ARTICLEAUTHORS{%
\AUTHOR{\bf Yuan Li}
\AFF{Georgia Institute of Technology, \EMAIL{yuanli@gatech.edu}}
\AUTHOR{\bf David A. Goldberg}
\AFF{Georgia Institute of Technology, \EMAIL{dgoldberg9@isye.gatech.edu}, \URL{http://www2.isye.gatech.edu/~dgoldberg9/}}
}
\RUNAUTHOR{{Yuan Li and David A. Goldberg}}


\ABSTRACT{%
We consider the FCFS $GI/GI/n$ queue in the Halfin-Whitt heavy traffic regime, in the presence of heavy-tailed distributions (i.e. infinite variance).  We prove that under minimal assumptions, i.e. only that processing times have finite $1 + \epsilon$ moment for some $\epsilon > 0$ and inter-arrival times have finite second moment, the sequence of stationary queue length distributions, normalized by $n^{\frac{1}{2}}$, is tight in the Halfin-Whitt regime.  All previous tightness results for the stationary queue length required that processing times have finite $2 + \epsilon$ moment.  Furthermore, we develop simple and explicit bounds on the stationary queue length in that setting.
\\\\When processing times have an asymptotically Pareto tail with index $\alpha \in (1,2)$, we bound the large deviations behavior of the limiting process (defined as any suitable subsequential limit), and derive a matching lower bound when inter-arrival times are Markovian.  Interestingly, we find that the large deviations behavior of the limit has a \emph{sub-exponential} decay, differing fundamentally from the exponentially decaying tails known to hold in the light-tailed setting, and answering an open question from \cite{GG13}.
\\\\For the setting where instead the inter-arrival times have an asymptotically Pareto tail with index $\alpha \in (1,2)$, we extend recent results of \cite{Reed.16} (who analyzed the case of deterministic processing times) by proving that for general processing time distributions, the sequence of stationary queue length distributions, normalized by $n^{\frac{1}{\alpha}}$, is tight (here we use the scaling of \cite{Reed.16}, which we refer to as the Halfin-Whitt-Reed scaling regime).  We are again able to bound the large-deviations behavior of the limit, and find that our derived bounds do not depend on the particular processing time distribution, and are in fact tight even for the case of deterministic processing times.
\\\\Our proofs proceed by extending the stochastic comparison approach of \cite{GG13}, and associated recent explicit bounds for multi-server queues formulated in \cite{G17a}, to the heavy-tailed setting.
}



\maketitle

%


\section{Introduction.}\label{Section:Introduction}
\subsection{Halfin-Whitt regime and literature review.}
The staffing of large-scale queueing systems, and the associated trade-offs, are a fundamental problem in Operations Research.  The insight that in many settings of interest one should scale the number of servers to exceed the arrival rate by a quantity on the order of the square-root of the arrival rate, i.e. the so-called square-root staffing rule, is by now well-known.  This setting is formalized by the so-called Halfin-Whitt scaling regime for parallel server queueing systems, studied originally by \citet{E.48} and \citet{J.74}, and formally introduced by \citet{HW.81}, who studied the $GI/M/n$ system (for large $n$) when the traffic intensity $\rho$ scales like $1 - B n^{-\frac{1}{2}}$ for some strictly positive excess parameter $B$.  There the authors prove weak convergence of the resulting queue-length process over compact time intervals, as well as weak convergence of the corresponding sequence of steady-state queue length distributions, when the queue-length of the $n$th system is normalized by $n^{\frac{1}{2}}$.  Namely, in both the transient and steady-state regimes, the queue-length scales like $n^{\frac{1}{2}}$ in the Halfin-Whitt regime when processing times are Markovian (and inter-arrival times have e.g. finite second moment).  We note that by queue-length, we refer to the number waiting in queue, not counting those jobs in service.
\\\indent The original results of \cite{HW.81} have since been extended in many directions.  Here we only review those results most relevant to our own investigations, and refer the interested reader to \cite{G16} for a comprehensive overview.  The most general results in the transient regime are those of \cite{R.09,PR.10}, which (customized to the setting of our own investigations, i.e. single-class parallel multi-server queues with i.i.d. inter-arrivals and processing times) prove that as long as the inter-arrival process satisfies a form of the central limit theorem on the scaling of $n^{\frac{1}{2}}$ (which will in general hold if the inter-arrival times have finite variance), and the processing time distribution has finite mean, then the associated sequence of queue-length processes, normalized by $n^{\frac{1}{2}}$, converges weakly to a non-trivial limiting process (if the system is initialized appropriately), described implicitly as the solution to a certain stochastic convolution equation.  
\\\indent As regards the scaling of the corresponding sequence of steady-state queue lengths, the most general known results are as follows.  Assuming that inter-arrival times and processing times have finite $2 + \epsilon$ moment for some $\epsilon > 0$, \cite{GG13} proves that the associated sequence of steady-state queue-lengths, normalized by $n^{\frac{1}{2}}$, is tight.  Under several additional technical assumptions, including that the processing times have finite third moment, the very recent results of \cite{AR15,AR16} show that the associated sequence has a unique weak limit.  Such a result was previously shown for the setting of processing times with finite support in \cite{GM08}.  In the presence of Markovian abandonments, an analogous result has been proven for the case of phase-type processing times.  Indeed, in this setting \cite{DDG14} proved that the sequence of steady-state queue length distributions, normalized by $n^{\frac{1}{2}}$, is tight with an explicit weak limit which the authors characterize as an Ornstein-Uhlenbeck process with piece-wise linear drift.  We note that although e.g. phase-type distributions are dense within the family of all distributions, due to the nature of the limits involved with the Halfin-Whitt regime, it is typically not clear how to translate results for such a restricted class of distributions to the general setting.
\\\indent Outside of the case of exponentially distributed processing times, the known characterizations for the limiting process (when such a limit is known to exist) are quite complicated.  As such, considerable effort has gone into understanding certain properties of this limit, where many of these results have pertained to the large deviations behavior of the limiting process.  In particular, for the case of inter-arrival times with finite second moment and processing times with finite support, \cite{GM08} prove that the the weak limit (associated with the sequence of normalized steady-state queue lengths) has an exponential tail, with a precise exponent identified as $-\frac{2B}{c^2_A + c^2_S}$, where $c^2_A (c^2_S)$ is the squared coefficient of variation (s.c.v) of inter-arrival (processing) times.  Namely, they prove that under those assumptions, the associated weak limit $\hat{Q}$ satisfies $\lim_{x \rightarrow \infty} x^{-1} \log\bigg( \pr\big( \hat{Q} > x \big) \bigg) = -\frac{2B}{c^2_A + c^2_S}$.  Put another way, the probability that the limiting process exceeds a large value $x$ behaves (roughly up to exponential order) like $\exp\big(-\frac{2B}{c^2_A + c^2_S} x\big)$.  The known results for the case of exponentially distributed and $H^*_2$ processing times yields the same exponent.  The stochastic comparison approach of \cite{GG13} was able to prove that the same exponent yields an upper bound on the large deviations behavior of any subsequential limit of the associated sequence of normalized queue-length random variables assuming only that inter-arrival and processing times have finite $2 + \epsilon$ moments for some $\epsilon > 0$, with equality for the case of exponentially distributed inter-arrival times.
\\\indent There has also been considerable interest in understanding the quality of Halfin-Whitt type approximations for finite $n$ (as opposed to having results which only hold asymptotically).  Such results include \cite{BD16,Gur14,BD15,Gur13,BDF15,G17a}.  We refer the reader to \cite{G17a} for a detailed overview and discussion, and note that none of these results apply to the heavy-tailed setting.  The very recent results of \cite{G17a} provided the first simple and explicit bounds for multi-server queues that scale universally as $\frac{1}{1 - \rho}$ across different notions of heavy traffic, including the Halfin-Whitt scaling.  However, the main results of \cite{G17a} assumed that both inter-arrival and processing times have finite $2 + \epsilon$ moment for some $\epsilon > 0$.
\subsection{Heavy tails in the Halfin-Whitt regime.}
A key insight from modern queueing theory is that when inter-arrival or processing times have a heavy tail (i.e. the tail of the probability distribution does not decay exponentially), the underlying system behaves qualitatively different, e.g. it may exhibit
long-range dependencies over time, and have a higher probability of rare events \cite{Gross.99}. As many applications in modern service systems (e.g. length of stay in a hospital, length of time of a call) are potentially highly variable (e.g. due to prolongued illnesses, or having to resolve a complex IT problem), and may experience traffic which is bursty in nature (e.g. long periods of low activity followed by periods of high activity) \cite{Bar.05}, and several studies have empirically verified this phenomena in applications relevant to the Halfin-Whitt scaling \cite{Brown.05,Mar.98}, it is important to understand how the presence of heavy tails changes the performance of parallel server queues in the Halfin-Whitt scaling regime.  Although there is a vast literature on parallel server queues with heavy-tailed inter-arrival and/or processing times (which we make no attempt to survey here, instead referring the reader to \cite{Res.07}), it seems that surprisingly, very little is known about how such systems behave qualitatively in the Halfin-Whitt regime.  
\\\indent We now survey what is known in this setting.  The results of \cite{R.09,PR.10} imply that when inter-arrival times have finite second moment (i.e. satisfy a classical central limit theorem) and processing times have finite mean (but may have infinite $1 + \epsilon$ moment for some $\epsilon \in (0,1)$), the associated sequence of transient queue-length processes, normalized by $n^{\frac{1}{2}}$, converges weakly (over compact time sets) to a non-trivial limiting process (if the system is initialized appropriately), described implicitly as the solution to a certain stochastic convolution equation. 
\\\indent \cite{Reed.16} considers the case in which inter-arrival times have (asymptotically) a so-called pure Pareto tail with index $\alpha \in (0,1)$, i.e. $\lim_{x \rightarrow \infty} x^{\alpha} \pr(A > x) = C$ for some $\alpha \in (1,2)$ and $C \in (0,\infty)$, and processing times are deterministic.  In this case, \cite{Reed.16} identifies a different scaling regime, a certain modification of the Halfin-Whitt scaling regime with the scaling modified to account for the heavy-tailed inter-arrivals.  In particular, Hurvich and Reed consider the associated sequence of $GI/D/n$ queues when the traffic intensity $\rho$ scales like $1 - B n^{-\frac{1}{\alpha}}$ for some strictly positive excess parameter $B$.  In this case, Reed proves that the associated sequence of steady-state waiting time random variables, rescaled so as to be multiplied by $n^{1 - \frac{1}{\alpha}}$, converges weakly to an explicit limiting distribution $\hat{W}$ characterized as the supremum of a certain infinite-horizon one-dimensional discrete-time random walk, i.e. a so-called $\alpha$-stable random walk, with drift $-B$.  Furthermore, although Reed does not explicitly prove it, it follows from an analysis nearly identical to that given in \cite{JMM.04} that by the distributional Little's law (which is applicable since processing times are deterministic), the sequence of steady-state queue-length distributions, normalized by $n^{\frac{1}{\alpha}}$, also converges to $\hat{W}$ (for completeness we will include a proof of this fact in our appendix).  Intuitively, the steady-state queue length in the $n$th system is thus approximately $\hat{W} n^{\frac{1}{\alpha}}$.  Namely, for $\alpha < 2$, $n^{\frac{1}{2}}$ \textbf{is no longer the correct scaling}.  This insight is quite interesting, although we note the important fact that Reed's results are restricted to the case of deterministic processing times.
\\\indent Essentially all other references in the literature to queues in the Halfin-Whitt regime with heavy tails are to open questions, which we now review.  In \cite{GG13}, the authors note that the identified limiting large deviations exponent $-\frac{2B}{c^2_A + c^2_S}$ equals zero when either inter-arrival or processing times have infinite variance, and leave as an open question identifying the correct behavior in the presence of heavy tails.  The question of tightness of the associated sequence of steady-state queue length distributions, normalized by $n^{\frac{1}{2}}$, is similarly left open when processing times have infinite variance.
\\\indent The very recent explicit bounds of \cite{G17a} for multi-server queues, which exhibit universal $\frac{1}{1-\rho}$ scaling across different heavy-traffic regimes (including the Halfin-Whitt scaling), left the extension to the heavy-tailed case open as well.  
However, we note that the results of \cite{BC98,BC99} prove that even for the single-server queue, $\frac{1}{1-\rho}$ is no longer the correct scaling as $\rho \uparrow 1$ when processing times are heavy-tailed, where the correct scaling instead involves a different function of $\rho$ depending on the tail of the processing time distribution.  Intriguingly, the transient results of  \cite{R.09,PR.10} show that in the Halfin-Whitt regime, even when processing times are heavy-tailed, $\frac{1}{1-\rho}$ is the correct scaling (at least for the transient queue-length distribution), as in the Halfin-Whitt regime $\frac{1}{1-\rho}$ will scale as the square root of the number of servers. As such, it seems that in the heavy-tailed setting, whether $\frac{1}{1-\rho}$ is the correct scaling depends heavily on precisely how one sends a sequence of queues into heavy traffic.  
\\\indent Indeed, it has been recognized in the literature that the order in which one takes limits plays a critical role in the heavy-tailed setting, e.g. when simultaneously looking at large-deviations behavior and heavy traffic, and such questions have been analyzed in \cite{Blanch11} for the single-server setting.  For the case of multiple servers, it is known that the interaction between the number of servers, the traffic intensity, and the large deviations behavior is very subtle \cite{Foss12}.  Recently, several results have been proven as regards the large deviations behavior when the number of servers and traffic intensity are held fixed \cite{Foss06,Foss12,BM15}.  However, much less is known as regards how the large deviations behavior scales while simultaneously altering the number of servers and traffic intensity.  Although some general explicit bounds are given in \cite{RSW13}, building on the earlier work of \cite{Scheller06}, those bounds do not scale properly in the Halfin-Whitt scaling, and e.g. depend sensitively on certain parameters being non-integer (with the bounds degrading as those parameters approach integers).  Several interesting bounds are also given in \cite{W.2000}, which proves that in certain settings heavy-tailed processing times lead to heavy-tailed waiting times.  However, the upper bounds presented there do not scale correctly in the Halfin-Whitt regime (see e.g. \cite{G17a} for a discussion of how bounds based on cyclic scheduling scale), while the implications of the proven lower bounds in the Halfin-Whitt regime are unclear.  We also note that using a robust-optimization approach, a different family of bounds was developed for a non-stochastic model of multi-server queues with heavy tails in \cite{Bandi15}, although those bounds also do not scale appropriately in the Halfin-Whitt regime.
\subsection{Questions for this work.}
The above discussion regarding heavy-tailed inter-arrival and processing times in the Halfin-Whitt regime motivates the following questions.
\begin{question}\label{Q1}
If the inter-arrival times have finite second moment but the processing times only have finite $1 + \epsilon$ moment for some $\epsilon \in (0,1)$, is the sequence of steady-state queue lengths in the Halfin-Whitt regime, normalized by $n^{\frac{1}{2}}$, still tight?  We note that a positive answer is known for the corresponding sequence of transient queue lengths (properly initialized) over a fixed compact time interval, but the corresponding question for the steady-state queues remains open.
\end{question}
\begin{question}\label{Q2}
Supposing that the answer to Question\ \ref{Q1} is yes, what can be said about the qualitative properties of the associated limiting process (technically any weak limit of the associated tight sequence of normalized steady-state queue lengths), e.g. what can be said about the large deviations behavior of such a limit?  This question becomes especially interesting in light of the large deviations exponent $- \frac{2B}{c^2_A + c^2_S}$ identified in all previous settings in the literature, which becomes zero in the case of infinite variance, and suggests that a fundamentally different behavior may arise.
\end{question}
\begin{question}\label{Q3}
For the setting in which inter-arrival times have infinite variance, can the scaling regime described by Reed in \cite{Reed.16}, henceforth referred to as the Halfin-Whitt-Reed scaling regime, be extended from the setting of deterministic processing times to the setting of general processing time distributions?  Do the same insights regarding tightness and asymptotic scaling hold?  Also, supposing the answer is yes, can anything be said about the qualitative properties, e.g. large deviations behavior, of the associated limits?
\end{question}
\begin{question}\label{Q4}
Is it possible to derive simple and explicit bounds for multi-server queues in the Halfin-Whitt regime, when processing times are heavy-tailed?  As all previous work on explicit, non-asymptotic bounds for queues in the Halfin-Whitt regime assumed processing times have a finite second moment, these would be the first such explicit bounds in the heavy-tailed setting.
\end{question}
\subsection{Our contribution.}
In this paper, we provide positive answers to Questions\ \ref{Q1} - \ref{Q4}.
\begin{answer}\label{Answer1}
We prove that, so long as inter-arrival times have finite second moment and processing times have finite $1 + \epsilon$ moment for some $\epsilon > 0$, the sequence of steady-state queue lengths, normalized by $n^{\frac{1}{2}}$, is tight.  Namely, the presence of heavy-tailed processing times does not interfere with the fact that the steady-state queue lengths scale like $n^{\frac{1}{2}}$ in the Halfin-Whitt regime.
\end{answer}
\begin{answer}\label{Answer2}
For the special case that the processing times have an asymptotically pure Pareto tail, i.e. $\lim_{x \rightarrow \infty} x^{\alpha} \pr(S > x) = C$ for some $\alpha \in (1,2)$ and $C \in (0,\infty)$, we explicitly bound the large deviations behavior of the corresponding limit.  In particular, we prove that the tail has a \emph{subexponential decay}, i.e. that $\limsup_{x \rightarrow \infty} x^{1 - \alpha} \log\bigg( \pr\big( \hat{Q} > x \big) \bigg)$ is at most an explicit strictly negative constant, for any weak limit $\hat{Q}$.  Furthermore, for the case of Markovian inter-arrival times, we prove a lower bound which certifies that this is indeed the exact large deviations behavior for any such weak limit.  Interestingly, in contrast to the light-tailed (i.e. finite variance) setting, here we find that rare events are fundamentally more likely, with the probability of seeing a large queue length $x n^{\frac{1}{2}}$ decaying like $\exp( - C' x^{\alpha - 1})$ with $\alpha - 1 \in (0,1)$ and $C'$ an explicit constant.  This in essence resolves the question of the previously identified large deviations exponent $- \frac{2B}{c^2_A + c^2_S}$ which vanishes in the infinite-variance setting, since $\lim_{x \rightarrow \infty} \frac{ C' x^{\alpha - 1} }{x} = 0$.  From a practical standpoint, this insight is important, as it suggests that when processing times are heavy-tailed (which as noted is a setting relevant in several service-system applications), it is much more likely to see large queue lengths, where we successfully quantify the meaning of ``much more likely".
\end{answer}
\begin{answer}\label{Answer3}
We prove that the Halfin-Whitt-Reed regime can indeed be extended to the setting of generally distributed processing times.  In particular, we prove that when inter-arrival times have an asymptotically pure Pareto tail with index $\alpha \in (1,2)$, and procesesing times have a finite $1 + \epsilon$ moment for some $\epsilon > 0$ (but are otherwise completely general), the sequence of steady-state queue lengths (under the Halfin-Whitt-Reed scaling), normalized by $n^{\frac{1}{\alpha}}$, is tight.  We also provide an explicit bound on the tail of the associated weak limit, and provide an upper bound on the associated large deviations exponent.  Intriguingly, we find that our upper bound closely resembles the exact weak limit proven for the special case of deterministic processing times by \cite{Reed.16}, and that in this special case the exact large deviations behavior actually matches our upper bound.
\end{answer}
\begin{answer}\label{Answer4}
We extend the framework of \cite{G17a} to provide the first simple and explicit bounds for multi-server queues that scale correctly in the Halfin-Whitt regime when processing times are heavy-tailed.
\end{answer}
\subsection{Outline of rest of paper.}
The rest of the paper proceeds as follows.  We state our main results in Section\ \ref{mainsec}.  We prove our explicit bounds for multi-server queues in the Halfin-Whitt regime when processing times may be heavy tailed in Section\ \ref{explicitsec}.  We prove our large deviations bounds for the setting that inter-arrival times have finite variance and processing times are heavy-tailed in Section\ \ref{ldsec}.  We extend the analysis of Reed from the special case of deterministic processing times to the case of general processing times, i.e. generalizing the notion of the Halfin-Whitt-Reed regime, in Section\ \ref{HWRsec}.  We provide a summary of our results and directions for future research in Section\ \ref{concsec}.

\section{Main results.}\label{mainsec}
In this section we formally state our main results.  
\subsection{Additional notations.}
As our main emphasis will be on queues in the Halfin-Whitt (-Reed) regime, we will customize our notations to this setting.  Let us fix  non-negative random variables A and S, with $\E[A] = \E[S] = 1$.  In general $A$ will not be the actual inter-arrival distribution to the queueing system of interest - instead a certain rescaling of $A$, with the rescaling depending on which results we are proving, will be the actual inter-arrival distribution (this is largely done as a notational simplification / convenience).  Here we note that by a simple rescaling argument, assuming both A and S have mean 1 is without loss of generality (w.l.o.g.).  Let ${\mathcal N}_o ({\mathcal A}_o)$ denote an ordinary renewal process with renewal distribution $S (A)$, and $N_o(t) \big( A_o(t) \big)$ the corresponding counting processes.  Let $\lbrace {\mathcal N}_i, i \geq 1 \rbrace \bigg( \lbrace {\mathcal N}_{o,i}, i \geq 1 \rbrace \bigg)$ denote a mutually independent collection of equilibrium (ordinary) renewal processes with renewal distribution $S$; ${\mathcal A}$ an independent equilibrium renewal process with renewal distribution $A$; and $\lbrace N_i(t), i \geq 1 \rbrace \bigg( \lbrace N_{o,i}(t), i \geq 1 \rbrace \bigg)$, $A(t)$ the corresponding counting processes.  Here we recall that an equilibrium renewal process (with renewal distribution $X$) is one in which the first renewal interval is distributed as the equilibrium distribution associated with $X$, i.e. letting $R(X)$ denote a r.v. such that $\pr( R(X) > y) = \frac{1}{\E[X]} \int_y^{\infty} \pr( X > z ) dz$, the first renewal interval is distributed as $R(X)$.  Also, let $\sigma_A (\sigma_S)$ denote $\big(Var[A]\big)^{\frac{1}{2}} \Big(\big(Var[S]\big)^{\frac{1}{2}}\Big)$, and $c_A (c_S)$ also denote $\sigma_A (\sigma_S)$ (here the standard deviation equals the coefficient of variation as the mean equals one).  Also, let $\lbrace A_i, i \geq 1 \rbrace$ $(\lbrace S_i, i \geq 1 \rbrace$) denote the sequence of inter-event times in ${\mathcal A}_o ({\mathcal N}_o)$.  Let us evaluate all empty summations to zero, and all empty products to unity; and as a convention take $\frac{1}{\infty} = 0$ and $\frac{1}{0} = \infty$.  For an event ${\mathcal E}$, let $I({\mathcal E})$ denote the corresponding indicator function.  Unless stated otherwise, all processes should be assumed right-continuous with left limits (r.c.l.l.), as is standard in the literature.  Also, for two r.v.s $X,Y$, let $X \sim Y$ denote equivalence in distribution.  For a real number $x$, let $x^+ \stackrel{\Delta}{=} \max(x,0)$, and $\sign(x)$ denote the sign of $x$, i.e. the function that evaluates to -1 for $x < 0$, 0 for $x = 0$, and 1 for $x > 0$.  In addition, for $t > 0$, let $\Gamma(t) \stackrel{\Delta}{=} \int_0^{\infty} x^{t - 1} \exp(-x) dx$ denote the well-known $\Gamma$-function.  By the so-called Euler reflection principle, the $\Gamma$-function is defined for negative (non-integer) values as follows: for $t > 0$ and non-integer, $\Gamma(-t) = \pi\big( - t \Gamma(t) \sin(\pi t) \big)^{-1}$.  We refer the interested reader to \cite{Abram.64} for further properties of this function, e.g. the useful fact that $\Gamma(t + 1) = t \times \Gamma(t)$ for all real $t$ (excluding negative integers).  Also, let $N$ denote a standard normal r.v.  

\subsubsection{Notation for queues in the Halfin-Whitt(-Reed) regime.}
For $B > 0, \alpha > 1$, and $n > B^{\frac{\alpha}{\alpha-1}}$, let $\lambda_{n,B,\alpha} \stackrel{\Delta}{=} n - B n^{\frac{1}{\alpha}}$, and ${\mathcal Q}^n_{A,S,B,\alpha}$ denote the FCFS $GI/GI/n$ queue with inter-arrival distribution $A \lambda^{-1}_{n,B,\alpha}$, and processing time distribution $S$.  If for any given initial condition, the total number of jobs in ${\mathcal Q}^n_{A,S,B,\alpha}$ (number in service + number waiting in queue) converges in distribution (as time goes to infinity, independent of the particular initial condition) to a steady-state r.v. $Q^n_{A,S,B,\alpha}(\infty)$, we say that ``$Q^n_{A,S,B,\alpha}(\infty)$ exists".  Here we refer the interested reader to \cite{AS08} for natural and mild technical conditions ensuring such existence.  For $n$ large, ${\mathcal Q}^n_{A,S,B,2}$ is said to be in the Half-Whitt (a.k.a. Quality-and-efficiency driven, QED) scaling regime \cite{HW.81}.  As Reed had studied ${\mathcal Q}^n_{A,S,B,\alpha}$ for $n$ large when $S$ is determinstic and $\alpha \in (1,2)$, we will generally say that ${\mathcal Q}^n_{A,S,B,\alpha}$ is in the Halfin-Whitt-Reed regime when $n$ is large and $\alpha \in (1,2)$.  In that case, supposing $Q^n_{A,S,B,\alpha}(\infty)$ exists, let us define $L^n_{A,S,B,\alpha}(\infty) \stackrel{\Delta}{=} \bigg( Q^n_{A,S,B,\alpha}(\infty) - n \bigg)^+$, i.e. the steady-state number of jobs waiting in queue (not counting those jobs in service).  Also, if for any given initial condition, the waiting time (i.e. time in system between time of arrival and time at which processing begins) for the $j$th job to arrive to $\mathcal{Q}^n_{A,S,B,\alpha}$ converges in distribution (as $j \rightarrow \infty$, independent of particular initial condition) to a steady-state r.v. $W^n_{A,S,B,\alpha}(\infty)$, we say that ``$W^n_{A,S,B,\alpha}(\infty)$ exists".
\subsection{Main results.}
\ \\\indent We begin by formalizing Answers\ \ref{Answer1} and\ \ref{Answer4}, i.e. stating our simple and explicit  bounds, as well as the implied tightness results.   We note that our tightness results are essentially the best possible, as the results of \cite{Reed.16} show that when inter-arrival times have infinite second moment square-root scaling is no longer appropriate.  In particular, the only case left unresolved is that in which $\E[S] < \infty$ but $\E[S^{1 + \epsilon}] = \infty$ for all $\epsilon > 0$.  Furthermore, even in that case, we believe our techniques could be extended to prove tightness and explicit bounds. 
\begin{theorem}[Answer\ \ref{Answer4}]\label{mainresult1}
Suppose that $\E[A^2] < \infty$, and $\E[S^{1 + \epsilon}] < \infty$ for some $\epsilon \in (0,1]$ (higher moments may or may not exist).  Then for all $B > 0$ and $n > 4 B^2$ such that $Q^n_{A,S,B,2}(\infty)$ exists, it holds that for all $x \geq 16$,
$\pr\bigg(  n^{-\frac{1}{2}} L^n_{A,S,B,2}(\infty) \geq x \bigg)$ is at most 
$$10^{100} \bigg( \epsilon \big(1 - \E[\exp(-S)]\big) \bigg)^{-7} (10 \E[S^{1 + \epsilon}])^{\frac{1}{\epsilon}} (1 + \sigma^2_A) (B^{-1} + B^{-2}) x^{- \frac{\epsilon}{1 + \epsilon}}.$$
\end{theorem}

\begin{corollary}[Answer\ \ref{Answer1}]\label{maintight}
Suppose that $\E[A^2] < \infty$, $\E[S^{1 + \epsilon}] < \infty$ for some $\epsilon \in (0,1]$ (higher moments may or may not exist), and for some $B > 0$, $Q^n_{A,S,B,2}(\infty)$ exists for all sufficiently large $n$.  Then $\big\lbrace n^{-\frac{1}{2}} L^n_{A,S,B,2}(\infty) , n > 4 B^2 \big\rbrace$ is tight.
\end{corollary}

We note that the tail decay rate demonstrated in Theorem\ \ref{mainresult1}, $x^{- \frac{\epsilon}{1 + \epsilon}}$, is likely not optimal.  As discussed in \cite{G17a}, the work of \cite{Scheller06} in fact suggests that as $n$ increases, the correct tail decay rate (and hence number of moments which are finite) scales with $n$ in a subtle manner, although how those moments scale (e.g. in the Halfin-Whitt regime) is unclear.  We leave the formulation of tighter uniform bounds in this setting as an interesting direction for future research.
\\\\\indent We next formalize Answer\ \ref{Answer2}, i.e. our large deviations results when processing times are asymptotically Pareto and inter-arrival times have finite second moment.  We begin by formulating a particular set of assumptions which we will need to state our results (which should be taken in addition to any assumptions posited to hold throughout the entire paper, e.g. $\E[A] = \E[S] = 1$).  
\begin{definition}[GH1 Assumptions]\ 
\begin{itemize}
\item $\E[A^2] < \infty$;
\item There exists $\alpha_S \in (1,2)$ and $C_S \in (0,\infty)$ s.t. $\lim_{x \rightarrow \infty} x^{\alpha_S} \pr(S > x) = C_S$;
\item $\limsup_{t \downarrow 0} t^{-1} \big( \pr ( S \leq t ) - \pr(S = 0) \big) < \infty$;
\item For each fixed $B > 0$, $Q^n_{A,S,B,2}(\infty)$ exists for all sufficiently large $n$.
\end{itemize}
\end{definition}
We note that the GH1 Assumptions are satisfied with appropriate parameters (for example) when $A$ has finite second moment and $S$ is a standard Pareto r.v. with tail index in $(1,2)$.
Let 
$$C_{B,S} \stackrel{\Delta}{=}  - C_S^{-1} B^{ 3 - \alpha_S} (\frac{\alpha_S - 1}{3 - \alpha_S})^{2 - \alpha_S} (2 - \alpha_S),$$
where we note that $C_{B,S} < 0$ under the GH1 Assumptions.  Then Answer\ \ref{Answer2} may be formalized as follows.
\begin{theorem}[Answer\ \ref{Answer2}]\label{heavyld1}
Under the GH1 Assumptions, 
$$
					\limsup_{x\rightarrow \infty}x^{-(\alpha_S - 1)} \log \Bigg( \limsup_{n\rightarrow \infty}\pr \bigg(n^{-\frac{1}{2}} L^n_{A,S,B,2}(\infty) > x\bigg)  \Bigg)\le C_{B, S}.
				$$
			
				If in addition $A$ is exponentially distributed, namely the system is $M/GI/n$, then
				\begin{eqnarray*}
					&\ & \liminf_{x \rightarrow \infty} x^{-(\alpha_S - 1)} \log \Bigg( \liminf_{n \rightarrow \infty} \pr\bigg( n^{-\frac{1}{2}} L^n_{A,S,B,2}(\infty) > x \bigg) \Bigg)  
					\\&=& \limsup_{x \rightarrow \infty} x^{-(\alpha_S - 1)} \log \Bigg( \limsup_{n \rightarrow \infty} \pr\bigg( n^{-\frac{1}{2}} L^n_{A,S,B,2}(\infty) > x \bigg) \Bigg)
					\ \ \ =\ \ \ C_{B,S}.
				\end{eqnarray*}
\end{theorem}
Roughly, Theorem\ \ref{heavyld1} implies that when processing times are asymptotically Pareto with power law decay parameter $\alpha_S \in (1,2)$, the probability of the queue exceeding a large queue length $x n^{\frac{1}{2}}$ decays roughly as 
$\exp\big( C_{B,S} x^{\alpha_S - 1} \big)$, which (since $\alpha_S - 1 \in (0, 1)$ and $C_{B,S} < 0$) decays sub-exponentially.  Namely, rare events are much more likely in this setting, as opposed to the light-tailed setting analyzed in \cite{GG13}, for which the decay was exponential.  Note that $|C_{B,S}|$ is increasing in B and decreasing in $C_S$, and hence in some sense seeing large queue lengths become ``less likely" as B increases (``more likely" as $C_S$ increases), which makes sense as when B is large the system is less loaded (when $C_S$ is large extreme procesing times are more likely), where we note that a similar monotonicity was observed in \cite{GG13} (with analogous quantities in the light-tailed setting).  Interestingly, the variability of the inter-arrival times does not appear in $C_{B,S}$, in contrast to the exponent identified in \cite{GG13} for the light-tailed setting.  This fact, combined with the tightness of our bound for the case of Markovian inter-arrival times, suggests that the tail behavior dictated by Theorem\ \ref{heavyld1} should in fact hold for any inter-arrival distribution with finite second moment, although a proof seems beyond the reach of current techniques.  We also note that our results could likely be extended to the setting of heavy-tailed processing times with more general tail behavior, in which case the analogous results would involve e.g. appropriate slowly-varying functions (cf. \cite{Whitt.02}), although we leave such an extension for future research.

Finally, we formalize Answer\ \ref{Answer3}, extending the Halfin-Whitt-Reed regime to generally distributed processing times.  First, we formalize the Halfin-Whitt-Reed scaling regime through an appropriate set of assumptions.
\begin{definition}[HWR-$\alpha$ Assumptions]\ 
\begin{itemize}
\item $\alpha \in (1,2)$;
\item There exists $C_A \in (0,\infty)$ s.t. $\lim_{x \rightarrow \infty} x^{\alpha} \pr(A > x) = C_A$;
\item There exists $\epsilon \in (0,1]$ s.t. $\E[S^{1 + \epsilon}] < \infty$;
\item For each fixed $B > 0$, $Q^n_{A,S,B,\alpha}(\infty)$ and $W^n_{A,S,B,\alpha}(\infty)$ exist for all sufficiently large $n$.
\end{itemize}
\end{definition}
Then our formalization of Answer\ \ref{Answer3} is as follows.   We begin by introducing some additional definitions and notations.
For $\alpha \in (1,2)$, let $$C_{\alpha} \stackrel{\Delta}{=} (1 - \alpha) \big( \Gamma(2 - \alpha) \cos(\frac{\pi}{2} \alpha) \big)^{-1},$$ where we note that $C_{\alpha} \in (0,\infty)$ for all $\alpha \in (1,2)$.  We next define the family of so-called $\alpha$-stable distributions, where we note that many different parametrizations appear for these variables throughout the literature, and our parametrization is consistent with that given in \cite{Whitt.02} and \cite{Sam.94}.  Given stability (i.e. index) parameter $\alpha \in (1,2)$, scale parameter $\sigma > 0$, skewness parameter $\beta \in [-1,1]$, and shift (i.e. location) parameter $\mu \in (-\infty,\infty)$, the corresponding $\alpha$-stable r.v. $S_{\alpha}(\sigma,\beta,\mu)$ is uniquely defined by its characteristic function (for all real $\theta$)
$$\E\big[\exp\big(i \theta S_{\alpha}(\sigma,\beta,\mu)\big)\big] = \exp\bigg(- (\sigma|\theta|)^{\alpha} \big( 1 - i \beta \sign(\theta) \tan(\frac{\pi}{2} \alpha) \big) + i \mu \theta \bigg).$$
Similarly, we define $\hat{S}_{\alpha,\beta}(t)_{t \geq 0}$ to be the corresponding (standardized) stochastic process known as a standardized $(\alpha,\beta)$-stable Levy motion \cite{Whitt.02,Sam.94}, where $\hat{S}_{\alpha,\beta}(0) = 0$, and for all $s,t \geq 0$,
$$\hat{S}_{\alpha,\beta}(s + t) - \hat{S}_{\alpha,\beta}(s) \sim t^{\frac{1}{\alpha}} S_{\alpha}(1,\beta,0).$$
We note that Levy motion is a Levy process (i.e. has stationary and independent increments), and has sample paths in the D-space (i.e. may have jumps), and we refer the interested reader to \cite{Whitt.02,Sam.94} for further details surrounding these processes, such as the fact that $\hat{S}_{\alpha,\beta}(t)_{t \geq 0}$ has the same distribution (on the process level) as $- \hat{S}_{\alpha,- \beta}(t)_{t \geq 0}$. 
\\\\Then our formalization of Answer\ \ref{Answer3} is as follows.

\begin{theorem}[Answer\ \ref{Answer3}]\label{beyondreedup1}
Under the HWR-$\alpha$ assumptions, $\big\lbrace n^{-\frac{1}{\alpha}} L^n_{A,S,B,\alpha}(\infty) , n \geq 1 \big\rbrace$ is tight.  Furthermore, for all $x > 0$,
\begin{equation}\label{beyondreedeq1}
\limsup_{n \rightarrow \infty} \pr\bigg(  n^{-\frac{1}{\alpha}} L^n_{A,S,B,\alpha}(\infty) > x \bigg) \leq
\pr\Bigg( \sup_{t \geq 0} \bigg( -(\frac{C_A}{C_{\alpha}})^{\frac{1}{\alpha}} \hat{S}_{\alpha,1}(t) - B t \bigg) > x \Bigg).
\end{equation}
\end{theorem}

Note that our bound does not depend on the particulars of the processing time distribution at all.  As $-\hat{S}_{\alpha,1}(t)_{t \geq 0}$ is a so-called spectrally negative Levy process (i.e. all jumps are negative), it is well-known that $\sup_{t \geq 0} \bigg( -(\frac{C_A}{C_{\alpha}})^{\frac{1}{\alpha}} \hat{S}_{\alpha,1}(t) - B t \bigg)$ follows a simple exponential distribution (cf. \cite{Bing75,Port.89}).  In particular, we have the following corollary, which follows immediately from Theorem\ \ref{beyondreedup1}, the results of \cite{Port.89} (which explicitly characterize the parameter of this exponential distribution), and some straightforward algebra.

\begin{corollary}\label{beyondreedupcor}
Under the HWR-$\alpha$ assumptions, for all $x > 0$,
\begin{equation}\label{beyondreedeq1b}
\limsup_{n \rightarrow \infty} \pr\bigg(  n^{-\frac{1}{\alpha}} L^n_{A,S,B,\alpha}(\infty) > x \bigg) \leq
\exp\big( - \big(\frac{B}{C_A \alpha \Gamma(-\alpha)}\big)^{\frac{1}{\alpha - 1}} x \big),
\end{equation}
where the right-hand-side of (\ref{beyondreedeq1}) equals the right-hand-side of (\ref{beyondreedeq1b}).
\end{corollary}

We note that $ - \big(\frac{B}{C_A \alpha \Gamma(-\alpha)}\big)^{\frac{1}{\alpha - 1}} < 0$.  Intriguingly, the explicit result of Reed for the special case of deterministic processing times yields a weak limit whose complimentary c.d.f. is nearly identical to the right-hand-side of (\ref{beyondreedeq1}), the only difference being that the supremum is taken over positive integer times, instead of all positive real times.  In particular, the following result follows almost immediately from the results of \cite{Reed.16} (Reed actually proved the analogous results for waiting times, and for completeness we include a formal proof translating those results to the setting of steady-state queue length in the appendix).

\begin{theorem}[\cite{Reed.16}]\label{reed1}
Suppose the HWR-$\alpha$ assumptions hold, and in addition $S$ is deterministic (i.e. the queueing system is a $GI/D/n$ queue).  Then there is a dense subset ${\mathcal S}$ of ${\mathcal R}^+$ s.t. for all $x \in {\mathcal S}$,
\begin{equation}\label{reed1eq}
\lim_{n \rightarrow \infty} \pr\bigg( n^{-\frac{1}{\alpha}} L^n_{A,S,B,\alpha}(\infty) > x \bigg) = 
\pr\Bigg( \sup_{k \geq 0} \bigg( -(\frac{C_A}{C_{\alpha}})^{\frac{1}{\alpha}} \hat{S}_{\alpha,1}(k) - B k \bigg) > x \Bigg).
\end{equation}
\end{theorem}

In light of Theorem\ \ref{reed1}, our upper bound (holding for general processing time distributions) is in some sense nearly tight 
even for the very special case of deterministic processing times.  Indeed, it is well-known that for a process with stationary and independent increments, there are straightforward ways to neatly bound the gap between the all-time supremum and the supremum over integer times (cf. \cite{MZ06,Will87}).  For example, such an analysis can be used to prove that the large deviations behavior of our upper bound is matched for the special case of deterministic processing times, i.e. both exhibit the same exponential rate of decay.  For completeness, we include a proof in the appendix.

\begin{corollary}\label{ldreed}
Under the HWR-$\alpha$ Assumptions, 
$$
					\limsup_{x\rightarrow \infty} x^{-1} \log \Bigg( \limsup_{n\rightarrow \infty}\pr \bigg(n^{-\frac{1}{\alpha}} L^n_{A,S,B,\alpha}(\infty) > x\bigg)  \Bigg)\le - \big(\frac{B}{C_A \alpha \Gamma(-\alpha)}\big)^{\frac{1}{\alpha - 1}}.
				$$
			
				If in addition $S$ is deterministic, namely the system is $GI/D/n$, then
				\begin{eqnarray*}
					&\ & \liminf_{x \rightarrow \infty} x^{-1} \log \Bigg( \liminf_{n \rightarrow \infty} \pr\bigg( n^{-\frac{1}{\alpha}} L^n_{A,S,B,\alpha}(\infty) > x \bigg) \Bigg)  
					\\&=& \limsup_{x \rightarrow \infty} x^{-1} \log \Bigg( \limsup_{n \rightarrow \infty} \pr\bigg( n^{-\frac{1}{\alpha}} L^n_{A,S,B,\alpha}(\infty) > x \bigg) \Bigg)
					\ \ \ =\ \ \ - \big(\frac{B}{C_A \alpha \Gamma(-\alpha)}\big)^{\frac{1}{\alpha - 1}}.
				\end{eqnarray*}
\end{corollary}

Whether $- \big(\frac{B}{C_A \alpha \Gamma(-\alpha)}\big)^{\frac{1}{\alpha - 1}}$ is the correct exponent for any given processing time distribution remains an interesting open question, although our results would certainly suggest that this should be the case.
\section{Explicit bounds and proof of Theorem\ \ref{mainresult1}.}\label{explicitsec}
In this section we prove Theorem\ \ref{mainresult1}, from which our tightness result Corollary\ \ref{maintight} will immediately follow.  We proceed by extending the framework of \cite{GG13,G17a} to the heavy-tailed setting.  We begin by reviewing several relevant results.
\subsection{Review of bounds from \cite{GG13}.}
In \cite{GG13}, the authors prove that the steady-state queue length of a $GI/GI/n$ can be bounded from above (in distribution) by the supremum of a relatively simple one-dimensional random walk.  We note that although to simplify notations the authors of \cite{GG13} imposed the restriction that $\pr(A = 0) = \pr(S = 0) = 0$ (to preclude having to deal with simultaneous events), this restriction is unnecessary and the proofs of \cite{GG13} can be trivially modified to accomodate this setting.  As such, we state the relevant stochastic-comparison result of \cite{GG13} here without that unnecessary assumption, albeit customized to our particular setting (i.e. in terms of the Halfin-Whitt(-Reed) regime).  

\begin{theorem}[\cite{GG13} Theorem 3]\label{ggold}
Suppose that $B > 0, \alpha > 1$, $n > B^{\frac{\alpha}{\alpha-1}}$, and $Q^n_{A,S,B,\alpha}(\infty)$ exists.  Then for all $x \geq 0$, 

\begin{equation}\label{GG1}
\pr\bigg(  n^{-\frac{1}{\alpha}} L^n_{A,S,B,\alpha}(\infty) \geq x \bigg) \leq \pr\Bigg( n^{-\frac{1}{\alpha}} \sup_{t \geq 0} \bigg( A(\lambda_{n,B,\alpha}t) - \sum_{i=1}^n N_i(t) \bigg) \geq x \Bigg).
\end{equation}
\end{theorem}

The authors also prove that the steady-state queue length can be lower-bounded by a different type of supremum, essentially dual to that given in \ref{GG1} (with the supremum and probability operators interchanged), when inter-arrival times are Markovian.  As we will later need these results for several proofs, we state them here.  Let $Z_{n,B,\alpha}$ be a Poisson r.v. with mean $\lambda_{n,B,\alpha}$.  
\begin{theorem}[\cite{GG13} Theorem 4]\label{ggold22}
Under the same assumptions as Theorem\ \ref{ggold}, supposing in addition that A is exponentially distributed, it holds that for all $x \geq 0$, 

\begin{equation}\label{GG2}
\pr\bigg(  n^{-\frac{1}{\alpha}} L^n_{A,S,B,\alpha}(\infty) \geq x \bigg) \geq \pr( Z_{n,B,\alpha} \geq n) \times \sup_{t \geq 0} \pr\bigg( n^{-\frac{1}{\alpha}} \big( A(\lambda_{n,B,\alpha}t) - \sum_{i=1}^n N_i(t) \big) \geq x \bigg).
\end{equation}
\end{theorem}

\subsection{Review of upper bounds from \cite{G17a}.}

In \cite{G17a}, the authors derive simple and explicit bounds for multi-server queues, which scale universally as $\frac{1}{1-\rho}$ across different heavy-traffic regimes, under the assumption that both inter-arrival and processing times have finite second moment.  As intermediate results, they also derived several lemmas which yield very general conditional bounds, which do not require the assumption of finite second moment.  These conditional results are of the form ``if certain quantities relating to the central moments of pooled renewal processes can be bounded by $\ldots$, then certain suprema appearing in the right-hand-side of (\ref{GG1}) can be bounded by $\ldots$".  The approach taken in \cite{G17a} to apply these conditional bounds did require $\E[S^2] < \infty$.  Here we take a different approach, which will allow us to utilize these same conditional bounds even in the heavy-tailed setting.  First, we remind the reader of several results from \cite{G17a}, including these general conditional bounds.

\begin{lemma}[\cite{G17a} Lemma 6]\label{alltimepooled1}
Suppose that for some fixed $n \geq 1, C_1, C_2 > 0; r_1 > s_1 > 1$; and $r_2 > 2$:
				\begin{enumerate}[(i)]
					\item For all $t \geq 1$,
$$\E\big[|\sum_{i=1}^{n} N_i(t) - n t|^{r_1}\big] \leq C_1 n^{\frac{r_1}{2}} t^{s_1}.$$
					\item For all $t \in [0,1]$,
$$\E\big[|\sum_{i=1}^{n} N_i(t) - n t|^{r_2}\big] \leq C_2 \max\big( n t , (n t)^{\frac{r_2}{2}} \big).$$
\end{enumerate}
Then for all $\nu > 0$ and $\lambda \geq 8$,
$$\pr\Bigg( \sup_{t \geq 0} \bigg( n t - \sum_{i=1}^{n} N_i(t) - \nu t  \bigg) \geq \lambda \Bigg)$$ 
is at most
$$
\bigg( \frac{100 (r_1 + r_2)^3}{(s_1 - 1)(r_1 - s_1)(r_2 - 2)}\bigg)^{r_1 + r_2 + 2}
\bigg(  C_1 n^{\frac{r_1}{2}} \nu^{-s_1} \lambda^{- (r_1 - s_1)} + C_2 n^{\frac{r_2}{2}} (\lambda \nu)^{- \frac{r_2}{2}} \bigg).
$$
\end{lemma}
\ \\\indent Second, we recall a useful bound from \cite{G17a} which will verify the conditions needed to apply Lemma\ \ref{alltimepooled1} for the case $t \in [0,1]$.  In that regime the fact that processing times are heavy-tailed does not lead to any pathologies, and thus we can simply use the results proven in \cite{G17a}.  Later we will develop new bounds to handle the $t \geq 1$ regime, where the heavy tails significantly change the analysis required to apply Lemma\ \ref{alltimepooled1}.
\begin{lemma}[\cite{G17a} Lemma 18]\label{binomial2}
For all $k \geq 1, p \geq 2, t \in [0,1],$ and $\theta > 0$, 
\begin{equation}\label{focus1}
\E\bigg[ \big| \sum_{i=1}^k N_i(t) - k t \big|^p \bigg] \leq \exp(\theta) \big(\frac{10^5 p^4 }{ 1 - \E[\exp(- \theta S)] }\big)^{p + 2} \max\big( kt , (k t)^{\frac{p}{2}} \big).
\end{equation}
\end{lemma}

\subsection{Novel bound for variance of pooled heavy-tailed renewal processes.}
In this section, we prove a novel simple, explicit, and non-asymptotic bound for the variance of a heavy-tailed equilibrium renewal process, i.e. $Var[N_1(t)]$.  We note that for the case $\E[S^2] < \infty$, both the renewal function (i.e. $\E[N_o(t)]$), and the variance of $N_1(t)$, are understood fairly precisely, with fairly tight bounds known (especially under further assumptions e.g. finite third moment, cf. \cite{Daley.78,Daley.80,Lorden.70,G16}).  The correct asymptotic scaling is also known in the heavy-tailed setting, under additional assumptions such as that $S$ is regularly varying, and/or belongs to the domain of attraction of an appropriate stable law (cf. \cite{BO.03,Geluk.97,Mohan.76,Teugels.68,GK.03}), and in some of our later large deviation results we will use certain of these precise asymptotics.  We also note that the literature contains certain non-explicit general results regarding the central moments of $N_1(t)$ under minimal moment conditions, showing e.g. that $E[S^{1 + \epsilon}] < \infty$ implies that $\E[|\N_{o,1}(t) - t|^{1 + \epsilon}]$ is asymptotically sublinear in $t$ (cf. \cite{Gut09}), although these results do not seem amenable to our analysis.  Here we provide a different result (which is, to our knowledge, new) under minimal assumptions on $S$.  The result builds on an elegant bounding argument of \cite{F64}, and a well-known explicit integral representation for $Var[N_1(t)]$ (cf. \cite{Daley.78,Daley.80,Lorden.70,Whitt.02}).

\begin{lemma}\label{boundvarpoolrenew}
Suppose that $\E[S^{1 + \epsilon}] < \infty$ for some $\epsilon \in (0,1]$.  Then for all $t \geq 0$, it holds that
$$Var[N_1(t)] \leq (4 \E[S^{1 + \epsilon}])^{\frac{1}{\epsilon}} \big( t + t^{1 + \frac{1}{1 + \epsilon}} \big).
$$
\end{lemma}

Our proof proceeds by first expressing $Var[N_1(t)]$ in terms of an integral involving the renewal function, and then using a result of \cite{F64} to bound the renewal function (and the aforementioned integral).
We begin by stating the desired integral representation.
\begin{lemma}[\cite{Daley.78,Daley.80,Lorden.70,Whitt.02}]\label{varintegral}
For all $t \geq 0$, it holds that
$$Var[N_1(t)] = 2 \int_0^t \bigg( \big( \E[N_o(s)] + 1 - s \big) - \frac{1}{2} \bigg) ds.$$
\end{lemma}
We next state the appropriate result from \cite{F64}, customized to our own setting.  In particular, the following lemma follows immediately from \cite{F64} Theorem 2, by taking the function $h$ defined there to be $h(x) = x^{1 + \epsilon}$.
\begin{lemma}[\cite{F64} Theorem 2]\label{Farrell1}
Suppose that $\E[S^{1 + \epsilon}] < \infty$ for some $\epsilon \in (0,1]$.  Then for all $t \geq 0$, it holds that
\begin{equation}\label{fareq1}
t - 1 \leq \E[N_o(t)] \leq t - 1 + (\E[S^{1 + \epsilon}])^{\frac{1}{1 + \epsilon}} (\E[N_o(t)] + 1)^{\frac{1}{1 + \epsilon}}.
\end{equation}
\end{lemma}
We note that Lemma\ \ref{Farrell1} does not directly provide an easily used bound for $\E[N_o(s)] + 1 - s$, as the right-hand-side of (\ref{fareq1}) is essentially a ``recursive bound" for $\E[N_o(s)]$, i.e. $\E[N_o(s)]$ is bounded in terms of a different function of $\E[N_o(s)]$.  We now show how to use Lemma\ \ref{Farrell1} to provide explicit bounds for $\E[N_o(s)] + 1 - s$.  
\begin{corollary}\label{Farrell2}
Under the same assumptions as Lemma\ \ref{Farrell1}, for all $t \geq 0$,
$$\E[N_o(t)] + 1 - t \leq (2 \E[S^{1 + \epsilon}])^{\frac{1}{\epsilon}} (1 + t^{\frac{1}{1 + \epsilon}}).$$
\end{corollary}
\proof{Proof}
Let us fix $t \geq 0$.  Letting $Y_t \stackrel{\Delta}{=} \E[N_o(t)] + 1 - t$, we conclude from Lemma\ \ref{Farrell1} that 
\begin{equation}\label{Farrell21}
0 \leq Y_t \leq (\E[S^{1 + \epsilon}])^{\frac{1}{1 + \epsilon}} (Y_t + t)^{\frac{1}{1 + \epsilon}}.
\end{equation}
If $Y_t = 0$, we are done.  Thus suppose $Y_t > 0$.  Then (\ref{Farrell21}) implies that 
\begin{equation}\label{Farrell22}
Y_t^{\frac{\epsilon}{1 + \epsilon}} \leq (\E[S^{1 + \epsilon}])^{\frac{1}{1 + \epsilon}} (1 + \frac{t}{Y_t})^{\frac{1}{1 + \epsilon}}.
\end{equation}
We first prove that $Y_t \leq \max\bigg(t, \big(2 \E[S^{1 + \epsilon}]\big)^{\frac{1}{\epsilon}}\bigg)$.  Indeed, suppose for contradiction that $Y_t > \max\bigg(t, \big(2 \E[S^{1 + \epsilon}]\big)^{\frac{1}{\epsilon}} \bigg)$.  Then (\ref{Farrell22}) implies that 
$$\big(2 \E[S^{1 + \epsilon}]\big)^{\frac{1}{1 + \epsilon}} < (\E[S^{1 + \epsilon}])^{\frac{1}{1 + \epsilon}} 2^{\frac{1}{1 + \epsilon}},$$
itself a contradiction, thus proving the desired statement, which itself implies that 
\begin{equation}\label{Farrell23}
0 < Y_t \leq \big(2 \E[S^{1 + \epsilon}]\big)^{\frac{1}{\epsilon}} + t.
\end{equation}
Plugging (\ref{Farrell23}) into the right-hand-side of (\ref{Farrell21}), applying the subadditivity of the function $f(x) = x^{\frac{1}{1 + \epsilon}}$ (which follows from concavity), and the fact that $\E[S^{1 + \epsilon}] \geq 1$ (by Jensen's inequality since $\E[S] = 1$), we find that
\begin{eqnarray*}
Y_t &\leq& (\E[S^{1 + \epsilon}])^{\frac{1}{1 + \epsilon}} \bigg(\big(2 \E[S^{1 + \epsilon}]\big)^{\frac{1}{\epsilon}} + 2 t\bigg)^{\frac{1}{1 + \epsilon}}
\\&\leq& (\E[S^{1 + \epsilon}])^{\frac{1}{1 + \epsilon}} \bigg(\big(2 \E[S^{1 + \epsilon}]\big)^{\frac{1}{\epsilon} \times \frac{1}{1 + \epsilon}} + (2 t)^{\frac{1}{1 + \epsilon}} \bigg)
\\&\leq& 2^{\frac{1}{\epsilon} \times \frac{1}{1 + \epsilon}} \times (\E[S^{1 + \epsilon}])^{\frac{1}{1 + \epsilon} \times (1 + \frac{1}{\epsilon})} \times (1 + t^{\frac{1}{1 + \epsilon}})
\\&\leq& (2 \E[S^{1 + \epsilon}])^{\frac{1}{\epsilon}} (1 + t^{\frac{1}{1 + \epsilon}}),
\end{eqnarray*}
completing the proof.
\endproof
With Lemma\ \ref{varintegral} and Corollary\ \ref{Farrell2} in hand, we now complete the proof of Lemma\ \ref{boundvarpoolrenew}.
\proof{Proof}[Proof of Lemma\ \ref{boundvarpoolrenew}]
It follows from Lemma\ \ref{varintegral} and Corollary\ \ref{Farrell2} that for all $t \geq 0$,
\begin{eqnarray*}
Var[N_1(t)] &\leq& 2  (2 \E[S^{1 + \epsilon}])^{\frac{1}{\epsilon}} \int_0^t \big( 1 + s^{\frac{1}{1 + \epsilon}}\big) ds 
\\&\leq& (4 \E[S^{1 + \epsilon}])^{\frac{1}{\epsilon}} \big(t + t^{1 + \frac{1}{1 + \epsilon}} \big),
\end{eqnarray*}
completing the proof.
\endproof
\subsection{Proof of Theorem\ \ref{mainresult1}.}
In this section we complete the proof of Theorem\ \ref{mainresult1}.  We begin by applying a straightforward union bound 
to the right-hand-side of (\ref{GG1}), along with non-negativity and some basic monotonicities, to conclude the following.
\begin{lemma}\label{union1}
Suppose that $B > 0, \alpha \in (1,2]$, and $n > B^{\frac{\alpha}{\alpha-1}}$.  Then for all $x \geq 0$,
$\pr\Bigg( n^{-\frac{1}{\alpha}} \sup_{t \geq 0} \bigg( A(\lambda_{n,B,\alpha}t) - \sum_{i=1}^n N_i(t) \bigg) \geq x \Bigg)$ is at most
\begin{eqnarray}
\ &\ &\ \pr\Bigg( n^{-\frac{1}{\alpha}} \sup_{t \geq 0}\bigg( A\big( \lambda_{n,B,\alpha} t \big) - (n - \frac{1}{2} B n^{\frac{1}{\alpha}}) t \bigg) \geq \frac{1}{2} x \Bigg) \label{reedeq1}
\\&\ &\ \ \ +\ \ \ \pr\Bigg( n^{-\frac{1}{2}} \sup_{t \geq 0}\bigg( \big( n t - \sum_{i=1}^n N_i(t) \big) - \frac{B}{2} n^{\frac{1}{2}} t \bigg) \geq \frac{1}{2} x \Bigg). \label{reedeq2}
\end{eqnarray}
\end{lemma}
\subsubsection{Bounding (\ref{reedeq1}), the supremum associated with the arrival process.}
In this section we bound (\ref{reedeq1}).  As here we want the most general result possible (i.e. only assuming finite second moment for the inter-arrival time distribution), we will proceed by relating the supremum to the waiting time in an appropriate single-server queue and applying Kingman's bound (as opposed to e.g. the analysis in \cite{G17a} which required stronger moment assumptions).  
We begin with a simple observation, following from the basic properties of ordinary and equilibrium renewal processes.  For $y > 1$, let $W_y$ denote a r.v. distributed as the steady-state waiting time in a $GI/GI/1$ queue with inter-arrival times distributed as $y A$ and processing times the constant 1.
\begin{observation}.\label{obsq1}
Suppose that $B > 0, \alpha \in (1,2]$, and $n > B^{\frac{\alpha}{\alpha-1}}$.  Then for all $\nu > \lambda_{n,B,\alpha}$ and $z \geq 0$, 
\begin{equation}\label{showsmall}
\pr\Bigg( n^{-\frac{1}{\alpha}} \sup_{t \geq 0}\bigg( A\big( \lambda_{n,B,\alpha} t \big) - \nu t \bigg) \geq z \Bigg)
\end{equation}
is at most
$$
\pr\Bigg( n^{-\frac{1}{\alpha}} \sup_{k \geq 0}\big( k - \nu \sum_{i=1}^k \frac{A_i}{\lambda_{n,B,\alpha}} \big) \geq z - n^{-\frac{1}{\alpha}} \Bigg).$$
It follows from Lindley's representation of the steady-state waiting time that (\ref{showsmall}) is at most 
$$
\pr\Bigg( n^{-\frac{1}{\alpha}} W_{\frac{\nu}{\lambda_{n,B,\alpha}}} \geq z - n^{-\frac{1}{\alpha}} \Bigg).$$
\end{observation}
Next, we recall the celebrated Kingman's bound for waiting times in a $GI/GI/1$ queue, only stating the result as customized to our particular setting.
\begin{lemma}[\cite{K62}, Kingman's Bound]\label{king1}
Suppose that $\E[A^2] < \infty$.  Then for all $y > 1$,
$$\E[W_y] \leq \frac{y^2 \sigma^2_A}{2 (y - 1)}.$$
\end{lemma}
Combining Observation\ \ref{obsq1} (with $\nu = n - \frac{1}{2} B n^{\frac{1}{\alpha}}$), Lemma\ \ref{king1} (with $y = \frac{n - \frac{1}{2} B n^{\frac{1}{\alpha}}}{n - B n^{\frac{1}{\alpha}}}$), Markov's inequality, and some straightforward algebra (e.g. the fact that $x \geq 4$ implies $\frac{x}{2} - n^{-\frac{1}{\alpha}} \geq \frac{x}{4}$, and $n > (2 B)^{\frac{\alpha}{\alpha-1}}$ implies $\frac{n - \frac{1}{2} B n^{\frac{1}{\alpha}}}{n - B n^{\frac{1}{\alpha}}} \leq 2$), we derive the following bound for (\ref{reedeq1}).
\begin{lemma}\label{boundapart}
Suppose that $\E[A^2] < \infty, B > 0, \alpha \in (1,2]$, and $n > (2 B)^{\frac{\alpha}{\alpha-1}}$.  Then for all $x \geq 4$, (\ref{reedeq1}) is at most 
$$10^2 \sigma^2_A B^{-1} n^{1 - \frac{2}{\alpha}} x^{-1}.$$
\end{lemma}
\subsubsection{Bounding (\ref{reedeq2}), the supremum associated with the departure process.}
We proceed by using Lemmas\ \ref{binomial2} and\ \ref{boundvarpoolrenew} to verify that the conditions of Lemma\ \ref{alltimepooled1} hold for appropriate parameters, which we use to bound (\ref{reedeq2}).  In particular, we prove the following.
\begin{lemma}\label{spartbound}
Suppose that $\E[S^{1 + \epsilon}] < \infty$ for some $\epsilon \in (0,1)$.  Then for all $B > 0, n \geq 1$, and $x \geq 16$, (\ref{reedeq2}) is at most
$$
10^{92} \epsilon^{-7} (8 \E[S^{1 + \epsilon}])^{\frac{1}{\epsilon}} \big(1 - \E[\exp(- S)] \big)^{-5} (B^{-1} + B^{-2}) x^{- \frac{\epsilon}{1 + \epsilon}}.$$
\end{lemma}
\proof{Proof}
By Lemma\ \ref{boundvarpoolrenew}, we find that for all $t \geq 1$,
$$\E\big[|\sum_{i=1}^{n} N_i(t) - n t|^{2}\big] \leq (8 \E[S^{1 + \epsilon}])^{\frac{1}{\epsilon}} n t^{1 + \frac{1}{1 + \epsilon}}.$$
By Lemma\ \ref{binomial2}, applied with $k = n, p = 3, \theta = 1$, we find that for all $t \in [0,1]$,
$$\E\big[|\sum_{i=1}^{n} N_i(t) - n t|^{3}\big] \leq \big(\frac{10^8 }{ 1 - \E[\exp(- S)] }\big)^{5} \max\big( n t , (n t)^{\frac{3}{2}} \big).$$
Thus we find that the conditions of Lemma\ \ref{alltimepooled1} are met with 
$$C_1 = (8 \E[S^{1 + \epsilon}])^{\frac{1}{\epsilon}}\ \ \ ,\ \ \ C_2 = \big(\frac{10^8 }{ 1 - \E[\exp(- S)] }\big)^{5}\ \ \ ,\ \ \ 
r_1 = 2\ \ \ ,\ \ \ s_1 = 1 + \frac{1}{1 + \epsilon}\ \ \ ,\ \ \ r_2 = 3.$$
Taking $\nu = \frac{B}{2} n^{\frac{1}{2}}$, $\lambda = \frac{x}{2} n^{\frac{1}{2}}$, we conclude that (\ref{reedeq2}) is at most 
\begin{eqnarray*}
\ &\ &\ \bigg( 10^6 \epsilon^{-1}\bigg)^{7} (8 \E[S^{1 + \epsilon}])^{\frac{1}{\epsilon}} \big(\frac{10^8 }{ 1 - \E[\exp(- S)] }\big)^{5}
\\&\ &\ \ \ \times\ \ \  \bigg(n\big( \frac{B}{2} n^{\frac{1}{2}} \big)^{-(1 + \frac{1}{1 + \epsilon})} \big( \frac{x}{2} n^{\frac{1}{2}} \big)^{-\frac{\epsilon}{1 + \epsilon}} + n^{\frac{3}{2}} \big(\frac{1}{4} x B n\big)^{-\frac{3}{2}} \bigg).
\end{eqnarray*}
Combining with some straightforward algebra completes the proof.
\endproof
With Lemma\ \ref{spartbound} in hand, we now complete the proof of Theorem\ \ref{mainresult1}.
\proof{Proof of Theorem\ \ref{mainresult1}}
Letting $\alpha = 2$, using Lemma\ \ref{boundapart} to bound (\ref{reedeq1}), and Lemma\ \ref{spartbound} to bound (\ref{reedeq2}), we conclude from Theorem\ \ref{ggold} and Lemma\ \ref{union1} (after some straightforward algebra) that for all $B > 0$, $n > 4 B^2$, and $x \geq 16$,
$\pr\big( n^{-\frac{1}{2}} L^n_{A,S,B,2}(\infty) \geq x \big)$ is at most
\begin{eqnarray*}
\ &\ &\ 
10^{92} \epsilon^{-7} (8 \E[S^{1 + \epsilon}])^{\frac{1}{\epsilon}} \big(1 - \E[\exp(- S)] \big)^{-5} (B^{-1} + B^{-2}) x^{- \frac{\epsilon}{1 + \epsilon}}
\\&\ &\ \ +\ \ 10^2 \sigma^2_A B^{-1} x^{-1}.
\end{eqnarray*}
Combining with some straightforward algebra, and Theorem\ \ref{ggold}, completes the proof.
\endproof

\section{Large deviations when $\E[A^2] < \infty$ and processing times are heavy-tailed, and proof of Theorem\ \ref{heavyld1}.}\label{ldsec}
In this section, we prove our large deviations results for the setting in which $\E[A^2] < \infty$ and $S$ is asymptotically Pareto with infinite variance, i.e. Theorem\ \ref{heavyld1}.  Our proof proceeds in a manner analogous to the large deviations results proven in \cite{GG13}.  In particular, we will use our tightness result to prove that our bound(s) for $L^n_{A,S,B,2}(\infty)$ behave like certain Gaussian processes in the H-W regime, where we note that (as in \cite{GG13}) some care will have to be taken as these bounds generally have the form of suprema over an infinite time horizon.  We will then use known results from the theory of Gaussian processes and heavy-tailed renewal processes to derive the appropriate large deviations behavior.  

\subsection{Preliminary weak convergence results.}
Before embarking on the proof of Theorem\ \ref{heavyld1}, we establish some preliminary weak convergence results to aid in our analysis.  For an excellent review of weak convergence, and the associated spaces (e.g. $D[0,T]$) and topologies/metrics (e.g. uniform, $J_1$, $M_1$), we refer the reader to \cite{Whitt.02}.  Recall that a Gaussian process on $\reals$ is a stochastic process $Z(t)_{t \geq 0}$ s.t. for any finite set of times $t_1,\ldots,t_k$, the vector $\big( Z(t_1),\ldots,Z(t_k) \big)$ has a Gaussian distribution.  A Gaussian process $Z(t)_{t \geq 0}$ is known to have its finite dimensional distributions uniquely determined by its mean function $\E[Z(t)]_{t \geq 0}$ and covariance function $\E[Z(s) Z(t)]_{s,t \geq 0}$, and refer the reader to \cite{Ad90}, and the references therein for details on existence, continuity, etc.  Let $\aleph(t)_{t \geq 0}$ denote the w.p.1 continuous Gaussian process s.t. $\E[\aleph(t)] = 0, \E[\aleph(s)\aleph(t)] = c^2_A \min(s,t)$, namely a driftless Brownian motion.  Then we may conclude the following from the well-known Functional Central Limit Theorem (FCLT) for renewal processes (see \cite{Whitt.02} Theorem 4.3.2 and Corollary 7.3.1).
\begin{theorem}[\cite{Whitt.02} Theorem 4.3.2 and Corollary 7.3.1]\label{renewalfclt}
Under the GH1 Assumptions, for any $T \in [0,\infty)$, the sequence of processes $\lbrace n^{-\frac{1}{2}} \big( A(\lambda_{n,B,2} t) - \lambda_{n,B,2} t \big)_{0 \leq t \leq T}, n \geq 1 \rbrace$ converges weakly to $\aleph(t)_{0 \leq t \leq T}$ in the space $D[0,T]$ under the $J_1$ topology.  
\end{theorem}
We now give a weak convergence result for $\sum_{i=1}^n N_i(t)$, which is stated in \cite{Whitt.02} (see Theorem 7.2.3) and formally proven in \cite{Whitt.85} (see Theorem 2). 
\begin{theorem}[\cite{Whitt.02} Theorem 7.2.3, \cite{Whitt.85} Theorem 2]\label{renewalfclt2}
There exists a w.p.1 continuous Gaussian process ${\mathcal D}(t)_{t \geq 0}$ s.t. $\E[{\mathcal D}(t)] = 0, \E[{\mathcal D}(s){\mathcal D}(t)] = \E[\big( N_1(s) - s \big) \big( N_1(t) - t \big)]$ for all $s,t \geq 0$.  Furthermore, under the GH1 Assumptions, for any $T \in [0,\infty)$, the sequence of processes $\lbrace n^{-\frac{1}{2}} \big( \sum_{i=1}^n N_i(t) - n t \big)_{0 \leq t \leq T}, n \geq 1 \rbrace$ converges weakly to ${\mathcal D}(t)_{0 \leq t \leq T}$ in the space $D[0,T]$ under the $J_1$ topology.
\end{theorem}
Let ${\mathcal Z}_{\infty}(t)_{t \geq 0}$ denote the Gaussian process s.t. ${\mathcal Z}_{\infty}(t) = \aleph(t) - {\mathcal D}(t)$ for all $t \geq 0$, and ${\mathcal Z}_{\infty,B}(t)_{t \geq 0}$ denote the Gaussian process s.t. ${\mathcal Z}_{\infty,B}(t) = \aleph(t) - {\mathcal D}(t) - B t$ for all $t \geq 0$.  Existence and continuity of both these processes follows from Theorems\ \ref{renewalfclt} and\ \ref{renewalfclt2}, which further imply the following (as similarly noted in \cite{GG13}).
\begin{corollary}\label{renewalfclt3}
Under the GH1 Assumptions, for any $T \in [0,\infty)$, the sequence of processes $\lbrace n^{-\frac{1}{2}} \big( A(\lambda_{n,B,2} t) - \sum_{i=1}^n N_i(t) \big)_{0 \leq t \leq T}, n \geq 1 \rbrace$ converges weakly to ${\mathcal Z}_{\infty,B}(t)_{0 \leq t \leq T}$ in the space $D[0,T]$ under the $J_1$ topology.
\end{corollary}

\subsection{Preliminary large deviations results.}
Next, we will need to establish some results from the theory of large deviations of Gaussian processes and their suprema.  We note that the relationship between the large deviations of suprema of Gaussian processes and the large deviations of queueing systems is well known, and there is a significant literature studying the large deviations of such processes (e.g. \cite{Dieker.05}).  
We will rely heavily on the following result, proven in \cite{Dieker.05} Proposition 1, describing the large deviation behavior of the supremum of certain Gaussian processes.  We note that a special case of the same result, customized to the light-tail setting, was also used in \cite{GG13}.  Before stating the result, let us recall the definition of a regularly varying function.
\begin{definition}[Regularly varying function]
A function $f: {\mathcal R}^+ \rightarrow {\mathcal R}^+$ is regularly varying with index $\gamma$ if for all $t > 0$,
$\lim_{x \rightarrow \infty} \frac{ f(t x) }{f(x)} = t^{\gamma}$.
\end{definition}
\ \\\indent We note that as the complimentary c.d.f.s of heavy-tailed distributions are typically regularly varying, the analysis of regularly varying functions is pervasive in the study of heavy-tailed phenomena, and we refer the interested reader to \cite{BGT89} for an excellent overview of the subject.  Then the aforementioned large deviations result is as follows.
\begin{lemma}[\cite{Dieker.05} Proposition 1] \label{diekerprop1}
Suppose $\mathcal{G}(t)_{t \geq 0}$ is a centered, continuous Gaussian process with stationary increments, satisfying the following conditions.
\begin{itemize}
\item The associated variance function $\E[\mathcal{G}^2(t)]$ is continuous (on ${\mathcal R}^+$) and regularly varying with index $2 H$ for some $0< H < 1$.  
\item There exists $\epsilon > 0$ s.t. $\lim_{t\downarrow 0} \E[\mathcal{G}^2(t)] |\log(t)|^{1+\epsilon} = 0$.
\end{itemize}
Then for all $\beta > H$ and $c > 0$,
$$
\lim_{x\rightarrow \infty} \bigg( \frac{ \E[\mathcal{G}^2(x^{\frac{1}{\beta}})] }{x^2} \log \pr \big( \sup_{t\ge 0} \mathcal{Z}(t) - c t^\beta \ge x \big) \bigg) = -\frac{1}{2} c^{\frac{2H}{\beta}} (\frac{H}{\beta-H})^{-\frac{2H}{\beta}} (\frac{\beta}{\beta - H})^2.
$$
\end{lemma}

We now use Lemma\ \ref{diekerprop1} to analyze the large deviations behavior of ${\mathcal Z}_{\infty,B}(t)_{t \geq 0}$, by proving that ${\mathcal Z}_{\infty}(t)_{t \geq 0}$ satisfies the conditions of Lemma\ \ref{diekerprop1} for an appropriate parameter of regular variation.  The proof relies on certain known results regarding the variance of heavy-tailed renewal processes (cf. \cite{GK.03}).  In particular, we recall a useful result regarding the variance of heavy-tailed renewal processes.  Such results have been proven under considerable generality (e.g. even when the first moment does not exist, and for asymptotic scaling beyond the second moment), although here we state the result customized to our own purposes and assumptions.
\begin{lemma}[\cite{GK.03}, Proposition 2]\label{heavyrenewvar}
Under the GH1 assumptions, 
$$
\lim_{t \rightarrow \infty} \frac{Var[N_1(t)]}{t^{3 - \alpha_S}} = 2 \big( (\alpha_S - 1) (2 - \alpha_S) (3 - \alpha_S) \big)^{-1} C_S.
$$
\end{lemma}
With Lemma\ \ref{heavyrenewvar} in hand, we now prove that ${\mathcal Z}_{\infty}(t)_{t \geq 0}$ satisfies the conditions of Lemma\ \ref{diekerprop1} for an appropriate parameter of regular variation, deferring all proofs to the appendix.
\begin{lemma}\label{ldlemma0}
Under the GH1 Assumptions, ${\mathcal Z}_{\infty}(t)_{t \geq 0}$ satisfies the conditions of Lemma\ \ref{diekerprop1}, where 
$\E[\mathcal{Z}^2(t)]$ is regularly varying with index $3 - \alpha_S$.
\end{lemma}
Finally, we combine Lemmas\ \ref{diekerprop1}\ -\ \ref{ldlemma0} to prove the desired large deviation results for ${\mathcal Z}_{\infty,B}(t)_{t \geq 0}$, again deferring the proof to the appendix.  

\begin{lemma}\label{ldlemma2}
Under the GH1 Assumptions,
\begin{equation}\label{underghld1}
\lim_{x\rightarrow\infty} \Bigg( x^{1 - \alpha_S} \log\bigg( \pr\big( \sup_{t\ge 0} {\mathcal Z}_{\infty,B}(t) \geq x \big) \bigg) \Bigg)
= C_{B,S}.
\end{equation}
\end{lemma}

Next, we state an additional large deviation-type result, which corresponds to the probability that ${\mathcal Z}_{\infty,B}$ exceeds a large value at the single time at which it is most likely to exceed that value (which will connect to an appropriate lower bound for multi-server queues).  The utility of considering such a quantity, in conjunction with the classical notion of large deviations considered in Lemma\ \ref{ldlemma2}, is well-known in the large-deviations literature, and we refer the interested reader to \cite{GG13} for further discussion.  We again defer all proofs to the appendix.

\begin{lemma}\label{ldlemma2bb}
Under the GH1 Assumptions,

\begin{equation}\label{underghld2}
\lim_{x\rightarrow\infty} \Bigg( x^{1 - \alpha_S} \log\bigg( \sup_{t\ge 0} \pr\big( {\mathcal Z}_{\infty,B}(t) \geq x \big) \bigg) \Bigg)
= C_{B,S}.
\end{equation}

\end{lemma}

\subsection{Weak convergence of the all-time supremum, and proof of Theorem\ \ref{heavyld1}.}
We now complete the proof of Theorem\ \ref{heavyld1} by combining our above large deviation results with a proof that 
$\lbrace  n^{-\frac{1}{2}} \sup_{t \geq 0}\big( A(\lambda_{n,B,2} t) - \sum_{i=1}^n N_i(t) \big), n \geq 1 \rbrace$ converges weakly to $\sup_{t \geq 0} {\mathcal Z}_{\infty,B}(t)$, along with our stochastic comparison bounds.  We note that such a result is not immediate, as the framework of weak convergence (of stochastic processes) generally deals only with compact time intervals, so extra care must be taken to handle such an infinite time horizon.  We note that closely related ideas were used in the proof of Lemma 7 and Theorem 2 in \cite{GG13}, although their proofs made use of the processing time distribution having finite second moment, and our result is stated in considerably greater generality.  Also, for a broader discussion of how the large deviations of the pre-limit connect to the large deviations of the limiting process, and the fundamental limits of such a line of inquiry, we refer the interested reader to \cite{GS12}.  Now, we prove the following general result, giving sufficient conditions for such an interchange to hold.  We defer all relevant proofs to the appendix.
\begin{lemma}\label{supinter1}
Suppose that $\lbrace Y_n(t)_{t \geq 0}, n \geq 1 \rbrace$ is a sequence of stochastic processes on $D[0,\infty)$ with stationary increments, and that ${\mathcal Y}_{\infty}(t)_{t \geq 0}$ is a fixed stochastic process (also with stationary increments, on $D[0,\infty)$).  Suppose also that:
\begin{enumerate}
\item $Y_n(0) = 0$ w.p.1 for all $n \geq 1$;
\item $\lbrace \sup_{t \geq 0} Y_n(t) , n \geq 1 \rbrace$ is tight; \label{ldlist1}
\item For all $M > 0$, $\lim_{t \rightarrow \infty} \pr\big( {\mathcal Y}_{\infty}(t) \geq - M \big) = 0$; \label{ldlist2}
\item For each fixed $T > 0$, $\lbrace \sup_{0 \leq t \leq T} Y_n(t) , n \geq 1 \rbrace$ converges weakly to $\sup_{0 \leq t \leq T}{\mathcal Y}_{\infty}(t)$. \label{ldlist3}
\end{enumerate}
Then $\lbrace \sup_{t \geq 0} Y_n(t), n \geq 1 \rbrace$ converges weakly to $\sup_{t \geq 0} {\mathcal Y}_{\infty}(t)$.
\end{lemma}

With the above results in hand, we now complete the proof of Theorem\ \ref{heavyld1}, noting that our proof proceeds similarly to the proof of the analogous large deviations result (which assumed $\E[S^2] < \infty$) in \cite{GG13}.

\proof{Proof of Theorem\ \ref{heavyld1}}
We begin by noting that under the GH1 assumptions $\lbrace n^{-\frac{1}{2}} \big( A(\lambda_{n,B,2} t) - \sum_{i=1}^n N_i(t) \big), n \geq 1 \rbrace$ satisfies the conditions of Lemma\ \ref{supinter1}, with limiting stochastic process ${\mathcal Z}_{\infty,B}$.  
Indeed, condition (\ref{ldlist1}.) follows immediately from our proof of Theorem\ \ref{mainresult1}.  Condition (\ref{ldlist2}.) follows from Lemma\ \ref{heavyrenewvar}, since that lemma (along with the definition of ${\mathcal Z}_{\infty,B}$) implies that $\limsup_{t \rightarrow \infty} \frac{Var[{\mathcal Z}_{\infty,B}(t)]}{t^{3 - \alpha_S}} < \infty$, which (combined with the strictly negative linear drift of ${\mathcal Z}_{\infty,B}$ and a straightforward argument involving the normal distribution which we omit) implies condition (\ref{ldlist2}.).  Finally, Condition (\ref{ldlist3}.) follows from Corollary\ \ref{renewalfclt3}, along with the continuity of the supremum map in the J1 topology, and the fact that convergence in J1 implies convergence of all co-ordinate projections corresponding to times $t$ such that w.p.1 the limit process has no jump exactly at time t (which will in this case be all $t \geq 0$) \cite{Whitt.02}.  It thus follows from Lemma\ \ref{supinter1} that $\lbrace n^{-\frac{1}{2}} \sup_{t \geq 0} \big( A(\lambda_{n,B,2} t) - \sum_{i=1}^n N_i(t) \big) , n \geq 1 \rbrace$ converges weakly to $\sup_{t \geq 0} {\mathcal Z}_{\infty,B}(t)$.  It follows (e.g. from the Portmanteau Theorem) that for all $x \geq 0$,
$$\limsup_{n \rightarrow \infty} \pr\bigg( n^{-\frac{1}{2}} \sup_{t \geq 0} \big( A(\lambda_{n,B,2} t) - \sum_{i=1}^n N_i(t) \big) \geq x \bigg) \leq \pr\big( \sup_{t \geq 0} {\mathcal Z}_{\infty,B}(t) \geq x \big).$$
The first part of Theorem\ \ref{heavyld1} (i.e. the upper bound) then follows by combining with Lemma\ \ref{ldlemma2bb} and the stochastic comparison result Theorem\ \ref{ggold}.  
\\\\We now prove the second part of Theorem\ \ref{heavyld1}, i.e. the lower bound.  Thus suppose $A$ is exponentially distributed.  Then it follows from Theorem\ \ref{ggold22} that for all $x \geq 0$ and $t \geq 0$, $\liminf_{n \rightarrow \infty} \pr\big(  n^{-\frac{1}{2}} L^n_{A,S,B,2}(\infty) > x \big)$ is at least 
$$\liminf_{n \rightarrow \infty} \pr \big( Z_{n,B,2} > n \big) \times \liminf_{n \rightarrow \infty} \pr\bigg( n^{-\frac{1}{2}} \big( A( \lambda_{n,B,2} t ) - \sum_{i=1}^n N_i(t) \big) > x \bigg),$$
which by the convergence of the Poisson to the normal, Corollary\ \ref{renewalfclt3}, and the Portmanteau Theorem is at least
$$\pr\big(N \geq B\big) \times \pr\big({\mathcal Z}_{\infty,B}(t) > x\big).$$
Taking the supremum over all $t \geq 0$, we conclude that 
\begin{equation}\label{provelb1}
\liminf_{n \rightarrow \infty} \pr\big(  n^{-\frac{1}{2}} Q^n_{A,S,B}(\infty) > x \big) \geq \pr\big(N \geq B\big) \times \sup_{t \geq 0} \pr\big({\mathcal Z}_{\infty,B}(t) > x\big).
\end{equation}
Combining with Lemma\ \ref{ldlemma2bb} and a straightforward limiting argument (the details of which we omit) then completes the proof.
\endproof

\section{The Halfin-Whitt-Reed regime, and proofs of Theorem\ \ref{beyondreedup1}.}\label{HWRsec}
In this section, we generalize the analysis of Reed from \cite{Reed.16} to the case of general processing times, and call the corresponding scaling regime the Halfin-Whitt-Reed regime.  First, we will need some additional background on so-called $\alpha$-stable processes and the generalized central limit theorem.
\subsection{The generalized central limit theorem.}
The celebrated central limit theorem describes the behavior of normalized partial sums of i.i.d. random variables which have finite variance, and proves that the sequence of normalized sums converges in distribution to a standard normal r.v.  In this section we review the generalization of these results to the setting in which the variance is infinite.  Here we only state a special case which will suffice for our purposes, e.g. only treating the case involving a pure pareto tail, only treating non-negative r.v.s, only treating the case $\alpha \in (1,2)$, etc. 

\begin{theorem}[Generalized CLT (\cite{Whitt.02}, Theorem 4.5.1)]\label{genclt1}
Suppose that $\lim_{x \rightarrow \infty} x^{\alpha} \pr(A > x) = C \in (0,\infty)$ for some $\alpha \in (1,2)$.  Then $\lbrace n^{-\frac{1}{\alpha}} \sum_{i=1}^n (A_i - 1) , n \geq 1 \rbrace$ converges in distribution to $(\frac{C}{C_{\alpha}})^{\frac{1}{\alpha}} S_{\alpha}(1,1,0)$, and we say that $A$ belongs to the normal domain of attraction of this limiting r.v.
\end{theorem}

There is also an analogous version of the functional central limit theorem for renewal processes.  

\begin{theorem}[Generalized FCLT for renewal processes (\cite{Whitt.02}, Corollary 7.3.2)]\label{genclt2}
Under the same assumptions as Theorem\ \ref{genclt1}, for any $T \in (0,\infty)$, $\big\lbrace n^{-\frac{1}{\alpha}}\big(A_o(n t) -  n t\big)_{0 \leq t \leq T} , n \geq 1 \big\rbrace$ and $\big\lbrace n^{-\frac{1}{\alpha}}\big(A(n t) -  n t\big)_{0 \leq t \leq T} , n \geq 1 \big\rbrace$ both converge weakly, in the space D[0,T] under the $M_1$ topology, to 
$$-(\frac{C}{C_{\alpha}})^{\frac{1}{\alpha}} \hat{S}_{\alpha,1}(t)_{0 \leq t \leq T}.$$
\end{theorem}

\subsection{Extending the Halfin-Whitt-Reed regime to general processing times, and proof of Theorem\ \ref{beyondreedup1}.}
In this section we use our stochastic-comparison approach, and results associated with our explicit bounds (i.e. Theorem\ \ref{mainresult1}), to extend the Halfin-Whitt-Reed regime beyond the case of deterministic process times.  In particular, we complete the proof of Theorem\ \ref{beyondreedup1}.  We proceed by means of a series of lemmas, and begin by proving the needed tightness result.
\begin{lemma}\label{reedtight}
Under the HWR-$\alpha$ assumptions, $\bigg\lbrace n^{-\frac{1}{\alpha}} \sup_{t \geq 0} \bigg( A(\lambda_{n,B,\alpha}t) - \sum_{i=1}^n N_i(t) \bigg) , n \geq 1 \bigg\rbrace$ is tight.
\end{lemma}
\proof{Proof}
By Lemma\ \ref{union1}, it suffices to prove tightness (separately) of 
\begin{equation}\label{showtighta}
\bigg\lbrace  n^{-\frac{1}{\alpha}} \sup_{t \geq 0}\bigg( A\big( \lambda_{n,B,\alpha} t \big) - (n - \frac{1}{2} B n^{\frac{1}{\alpha}}) t \bigg) \bigg\rbrace,
\end{equation}
and
\begin{equation}\label{showtightb}
\bigg\lbrace  n^{-\frac{1}{2}} \sup_{t \geq 0}\bigg( \big( n t - \sum_{i=1}^n N_i(t) \big) - \frac{1}{2} B n^{\frac{1}{2}} t \bigg) \bigg\rbrace.
\end{equation}
As tightness of (\ref{showtightb}) follows immediately from Lemma\ \ref{spartbound}, it suffices to demonstrate tightness of (\ref{showtighta}).  However, tightness of (\ref{showtighta}) follows immediately from Observation\ \ref{obsq1}, and Theorem 7.1 of \cite{BC99}, which gives sufficient conditions for tightness of the sequence of waiting times associated with a sequence of single-server queues with heavy-tailed inter-arrival times in heavy traffic.
\endproof
Next, we prove the appropriate weak convergence result.
\begin{lemma}\label{reedweak}
Under the HWR-$\alpha$ assumptions, for all $T \in (0,\infty)$,
$\bigg\lbrace n^{-\frac{1}{\alpha}} \sup_{t \in [0,T]} \bigg( A(\lambda_{n,B,\alpha}t) - \sum_{i=1}^n N_i(t) \bigg) , n \geq 1 \bigg\rbrace$ converges weakly to $\sup_{t \in [0,T]} \bigg( -(\frac{C}{C_{\alpha}})^{\frac{1}{\alpha}} \hat{S}_{\alpha,1}(t) - B t \bigg).$
\end{lemma}
\proof{Proof}
Note that $n^{-\frac{1}{\alpha}}\bigg( A(\lambda_{n,B,\alpha}t) - \sum_{i=1}^n N_i(t) \bigg)$ equals
\begin{eqnarray*}
\ &\ &\ (\frac{\lambda_{n,B,\alpha}}{n})^{\frac{1}{\alpha}} \times \lambda_{n,B,\alpha}^{-\frac{1}{\alpha}} \bigg( A( \lambda_{n,B,\alpha} t ) - \lambda_{n,b,\alpha} t \bigg)
\\\ &\ &\ \ \ +\ \ \ n^{\frac{1}{2} - \frac{1}{\alpha}} \times n^{-\frac{1}{2}} \big( n t - \sum_{i=1}^n N_i(t) \big)
\\\ &\ &\ \ \ -\ \ \ B t.
\end{eqnarray*}
Combining with Theorem\ \ref{genclt2},  Lemma\ \ref{spartbound}, and the basic properties of $J_1$ and $M_1$ convergence, e.g. continuity of the supremum map (cf. \cite{Whitt.02}) and well-known conditions for convergence of co-ordinate projections, along with the basic properties of spectrally negative $\alpha$-stable Levy processes, completes the proof of the desired weak convergence. 
\endproof
With Lemmas\ \ref{reedtight} and \ref{reedweak} in hand, we now complete the proof of Theorem\ \ref{beyondreedup1}.
\proof{Proof of Theorem\ \ref{beyondreedup1}}
In light of Theorem\ \ref{ggold}, it suffices to verify that 
$\bigg\lbrace n^{-\frac{1}{\alpha}} \bigg( A(\lambda_{n,B,\alpha}t) - \sum_{i=1}^n N_i(t) \bigg)_{t \geq 0} , n \geq 1 \bigg\rbrace$ satisfies the conditions of Lemma\ \ref{supinter1}.  In light of Lemmas\ \ref{reedtight} and\ \ref{reedweak}, it suffices to verify that for all $M > 0$, 
\begin{equation}\label{showmeee1}
\lim_{t \rightarrow \infty} \pr\bigg( -(\frac{C}{C_{\alpha}})^{\frac{1}{\alpha}} \hat{S}_{\alpha,1}(t) - B t \geq - M \bigg) = 0;
\end{equation}
equivalently (by the basic properties of $\alpha$-stable Levy processes) that for all $M > 0$,
$$\lim_{t \rightarrow \infty} \pr\bigg( - (\frac{C}{C_{\alpha}})^{\frac{1}{\alpha}} t^{\frac{1}{\alpha}} S_{\alpha}(1,1,0) \geq B t - M \bigg) = 0.$$
(\ref{showmeee1}) then follows from the fact that $\alpha > 1$, and $S_{\alpha}(1,1,0)$ is a.s. finite.  Combining the above verifies that the conditions of Lemma\ \ref{supinter1} are met, completing the proof.
\endproof

\section{Conclusion.}\label{concsec}
In this paper, we provided the first analysis of steady-state multi-server queues in the Halfin-Whitt regime when processing times have infinite variance.  We proved that under minimal assumptions, i.e. only that processing times have finite $1 + \epsilon$ moment for some $\epsilon > 0$ and inter-arrival times have finite second moment, the sequence of stationary queue length distributions, normalized by $n^{\frac{1}{2}}$, is tight in the Halfin-Whitt regime.  This confirmed that the presence of heavy tails in the processing time distributions does not change the fundamental scaling of the steady-state queue length, as $n^{\frac{1}{2}}$ was also the correct scaling in the light-tailed case, and was known to be the correct scaling for the transient queue length in the presence of heavy tails.  Furthermore, we developed simple, explicit, and uniform bounds for the steady-state queue length in the Halfin-Whitt regime, under only these minimal assumptions.  
\\\indent When processing times have an asymptotically Pareto tail with index $\alpha \in (1,2)$, we were able to bound the large deviations behavior of the limiting process (defined as any suitable subsequential limit), and derived a matching lower bound when inter-arrival times are Markovian.  Interestingly, we find that the large deviations behavior of the limit has a sub-exponential decay, differing fundamentally from the exponentially decaying tails known to hold in the light-tailed setting.  Also, for the setting where instead the inter-arrival times have an asymptotically Pareto tail with index $\alpha \in (1,2)$, we extended recent results of \cite{Reed.16} (who analyzed the case of deterministic processing times) by proving that for general processing time distributions, the sequence of stationary queue length distributions, normalized by $n^{\frac{1}{\alpha}}$, is tight (here we used the scaling of \cite{Reed.16}, which we named the Halfin-Whitt-Reed scaling regime).  Interestingly, our derived bounds do not depend at all on the specifics of the processing time distribution, and are nearly tight even for the case of deterministic processing times.  We further formalized this by using our results to prove a universal bound on the large deviations behavior of the associated limiting process, and proved that even the setting of deterministic processing times yields a matching large deviations exponent.
\\\indent Our work leaves several interesting directions for future research.  Within the Halfin-Whitt regime, there is the obvious question of deriving tighter explicit bounds, e.g. doing away with the massive constant appearing in our bounds, and developing tighter bounds on the demonstrated tail decay rate.  One could also extend our analysis to more general heavy-tailed distributions, e.g. not having an asymptotically pure Pareto tail, as well as analyze different queueing quantities (e.g. the steady-state probability of delay, for which some interesting results are proven in \cite{Reed.16}) to gain further insight into the impact of heavy tails on queues in the Halfin-Whitt(-Reed) regime.  Developing a deeper understanding of the weak limit process arising in the Halfin-Whitt(-Reed) regime (both with and without heavy tails), as well as how the processing time distribution impacts the behavior of this process, both remain largely open questions.  
\\\indent Even more interesting is the question of deriving any kind of simple and explicit bounds that scale universally across different notions of heavy traffic in the heavy-tailed setting, as was accomplished under the assumption of a finite $2 + \epsilon$ moment in \cite{G17a}.  On a related note, it would be very interesting to use our stochastic comparison approach to analyze the large deviations behavior of multi-server queues with heavy-tailed processing times for a fixed number of servers, where it is known that the interaction between the number of servers, the traffic intensity, and the large deviations behavior can be very subtle \cite{Foss12}.  Another question along these lines is to develop a clearer understanding of the connection (under e.g. the Halfin-Whitt scaling) between the finiteness of moments of the steady-state queue length, and how those moments scale with the traffic intensity.  Although the question of which moments are finite is by now fairly well understood \cite{Scheller06}, the question of how those finite moments scale in heavy traffic remains largely open, where we note that some interesting progress there follows from the recent results of \cite{G17a}.  
\\\indent More generally, developing a broad understanding of the connection between heavy tails, heavy traffic, large deviations, and e.g. the relative scaling of various quantities of interest remains an interesting open question for multi-server systems, especially if the number of servers is allowed to diverge as the traffic intensity approaches unity.  The same goes for our understanding of so-called sample-path large deviations, i.e. the question of the most likely way for such rare events to occur.  Indeed, at this time our understanding of such questions in the single-server setting (cf. \cite{Blanch11}) far outpaces our understanding in the multi-server case, which remains an interesting direction for future research.  
\\\indent On a final note, there is the important question of to what extent genuinely heavy-tailed phenomena arise in practice (e.g. in service systems), and how the resulting phenomena observed in practice connect to our theoretical understanding of heavy tails.  Answering such questions will no doubt require interdisciplinary work at the interface of statistics, probability, and (more broadly) data science, in the spirit of \cite{Brown.05}.

\section*{Acknowledgements.}
The authors gratefully acknowledge support from NSF grant no. 1333457, as well as several stimulating conversations with Ton Dieker, Kavita Ramanan, Josh Reed, and Bert Zwart.

\section{Appendix.}\label{appsec}
\subsection{Proof of Lemma\ \ref{ldlemma0}}

\proof{Proof of Lemma\ \ref{ldlemma0}.}
That $\mathcal{Z}(t)_{t\ge 0}$ is (w.p.1) continuous, centered, and has the stationary increments property follows from the corresponding properties of $\aleph(t)_{t \geq 0}$ and ${\mathcal D}(t)_{t \geq 0}$.  Since
\begin{equation}\label{varrep1}
\E[{\mathcal Z}^2(t)] = c^2_A t + Var[N_1(t)],
\end{equation}
continuity of $\E[{\mathcal Z}^2(t)]$, as well as the fact that $\lim_{t \downarrow 0} \E[{\mathcal Z}^2(t)] \log^2(t) = 0$,
follows from the integral representation Lemma\ \ref{varintegral}.  Combining with the regular variation implied by Lemma\ \ref{heavyrenewvar} completes the proof.
\endproof

\subsection{Proof of Lemma\ \ref{ldlemma2}.}
\proof{Proof of Lemma\ \ref{ldlemma2}}
It follows from Lemma\ \ref{ldlemma0} that under the GH1 Assumptions, we may apply Lemma\ \ref{diekerprop1} to $\sup_{t \geq 0} {\mathcal Z}_{\infty,B}(t)$, with ${\mathcal G}(t)_{t \geq 0} = {\mathcal Z}_{\infty}(t)_{t \geq 0}, c = B, \beta = 1, H = \frac{1}{2}(3 - \alpha_S)$.  It follows from Lemma\ \ref{heavyrenewvar} and (\ref{varrep1}) that (in the language of Lemma\ \ref{diekerprop1})
\begin{equation}\label{ldlemma2a1}
\lim_{x \rightarrow \infty} \bigg( \big(\frac{ \E[\mathcal{G}^2(x^{\frac{1}{\beta}})] }{x^2}\big) x^{\alpha_S - 1} \bigg) = 
2 \big( (\alpha_S - 1) (2 - \alpha_S) (3 - \alpha_S) \big)^{-1} C_S,
\end{equation}
and
\begin{equation}\label{ldlemma2a2}
-\frac{1}{2} c^{\frac{2H}{\beta}} (\frac{H}{\beta-H})^{-\frac{2H}{\beta}} (\frac{\beta}{\beta - H})^2
=
- 2 B^{3 - \alpha_S} (3 - \alpha_S)^{-(3 - \alpha_S)} (\alpha_S - 1)^{-(\alpha_S - 1)}.
\end{equation}
Combining with Lemma\ \ref{diekerprop1} and some straightforward algebra completes the proof.
\endproof

\subsection{Proof of Lemma\ \ref{ldlemma2bb}.}
\proof{Proof of Lemma\ \ref{ldlemma2bb}}
For $x \in {\mathcal R}^+$, let $T_{S,x} \stackrel{\Delta}{=} \frac{(3-\alpha_S) x}{B (\alpha_S-1)}$.  Note that 
\begin{eqnarray}
\ &\ &\ \liminf_{x\rightarrow\infty} \Bigg( x^{1 - \alpha_S} \log\bigg( \sup_{t\ge 0} \pr\big( {\mathcal Z}_{\infty,B}(t) \geq x \big) \bigg) \Bigg) \nonumber
\\&\ &\ \ \ \geq\ \ \ \liminf_{x\rightarrow\infty} \Bigg( x^{1 - \alpha_S} \log\bigg( \pr\big( {\mathcal Z}_{\infty,B}(T_{S,x}) \geq x \big) \bigg) \Bigg) \nonumber
\\&\ &\ \ \ =\ \ \ \liminf_{x\rightarrow\infty} x^{1 - \alpha_S} \log\Bigg( \pr\bigg( N >  2 (\alpha_S-1)^{-1} x \big(\E[{\mathcal Z}^2_{\infty,B}(T_{S,x})]\big)^{-\frac{1}{2}} \bigg) \Bigg). \label{underghld2eq1}
\end{eqnarray}
As it follows from Lemma\ \ref{ldlemma0} that $\lim_{x \rightarrow \infty} x \big(\E[{\mathcal Z}^2_{\infty,B}(T_{S,x})]\big)^{-\frac{1}{2}} = \infty$, and standard bounds for the normal distribution c.d.f. (cf. \cite{G16} Lemma 6) imply that there exists $y_0$ s.t. $y > y_0$ implies $\pr (N>y) \ge \exp (-\frac{y^2}{2} - y)$, we may further conclude that (\ref{underghld2eq1}) is at least
$$- \liminf_{x\rightarrow\infty} x^{1 - \alpha_S} \bigg(2 (\alpha_S - 1)^{-2} x^2 \big(\E[{\mathcal Z}^2_{\infty,B}(T_{S,x})]\big)^{-1} + 
2 (\alpha_S - 1)^{-1} x \big(\E[{\mathcal Z}^2_{\infty,B}(T_{S,x})]\big)^{-\frac{1}{2}} \bigg),$$
which by Lemma\ \ref{heavyrenewvar}, (\ref{varrep1}), and some straightforward algebra equals $C_{B,S}$.  Combining with the fact that, by the basic properties of the supremum operator, 
$$
\limsup_{x\rightarrow\infty} \Bigg( x^{1 - \alpha_S} \log\bigg( \sup_{t\ge 0} \pr\big( {\mathcal Z}_{\infty,B}(t) \geq x \big) \bigg) \Bigg)$$
is bounded (from above) by the left-hand-side of (\ref{underghld1}) completes the proof.
\endproof

\subsection{Proof of Lemma\ \ref{supinter1}.}
\proof{Proof of Lemma\ \ref{supinter1}}
First, we claim that
\begin{equation}\label{ldproofeq1}
\lim_{T \rightarrow \infty} \limsup_{n \rightarrow \infty} \pr \big( \sup_{ t \geq T } Y_n(t) \geq 0 \big) = 0.
\end{equation}
Indeed, for all $T > 0$, $M > 0$, and $n \geq 1$, by a union bound and stationary increments $\pr \big( \sup_{ t \geq T } Y_n(t) \geq 0 \big)$ is at most
\begin{equation}\label{unionsup1}
\pr\big( Y_n(T) \geq - M \big) + \pr\big( \sup_{t \geq 0} Y_n(t) \geq M \big).
\end{equation}
It follows from (\ref{ldlist1}.) - (\ref{ldlist2}.) that for any given $\epsilon > 0$, we may select $M_{\epsilon}, T_{\epsilon} \in (0,\infty)$ s.t. $\pr\big( {\mathcal Y}_{\infty}(T_{\epsilon}) \geq - M_{\epsilon} \big) < \frac{\epsilon}{2}$, and $\limsup_{n \rightarrow \infty} \pr\big( \sup_{t \geq 0} Y_n(t) \geq M_{\epsilon} \big) < \frac{\epsilon}{2}$.  Combining with (\ref{ldlist3}.), (\ref{unionsup1}), and the monotonicity of the supremum operator, it follows that for all $T \geq T_{\epsilon}$, $\limsup_{n \rightarrow \infty} \pr \big( \sup_{ t \geq T } Y_n(t) \geq 0 \big) < \epsilon$.  Combining with the definition of limit completes the proof of (\ref{ldproofeq1}).
\\\\It follows from (\ref{ldproofeq1}) that for any $x \geq 0$,
we may construct a strictly increasing sequence of integers $\lbrace T_{x,k^{-1}}, k \geq 1 \rbrace$ s.t. for all $k \geq 1$,
$$
\limsup_{n \rightarrow \infty} \pr \bigg( \sup_{ t \geq T_{x,k^{-1}} } Y_n(t) \geq x \bigg) < k^{-1}.
$$
Thus by a union bound, for all $x \geq 0$ and $k \geq 1$,
$$\limsup_{n \rightarrow \infty} \pr \bigg( \sup_{ t \geq 0 } Y_n(t) \geq x \bigg) \leq \limsup_{n \rightarrow \infty} \pr \bigg( \sup_{ 0 \leq t \leq T_{x,k^{-1}} } Y_n(t) \geq x \bigg) + k^{-1}.
$$
By letting $k \rightarrow \infty$, and applying (\ref{ldlist3}.), the monotonicity of the supremum operator, and the Portmanteau Theorem, we conclude that for all $x \geq 0$, 
\begin{equation}\label{firstparta}
\limsup_{n \rightarrow \infty} \pr \bigg( \sup_{ t \geq 0 } Y_n(t) \geq x \bigg) \leq \pr\big(\sup_{t \geq 0} {\mathcal Y}_{\infty}(t) \geq x \big).
\end{equation}
Next, we prove the analagous result for $\liminf_{n \rightarrow \infty} \pr \bigg( \sup_{ t \geq 0 } Y_n(t) > x \bigg)$.  In particular, for any fixed $T$, (\ref{ldlist3}.), the monotonicity of the supremum operator, and the Portmanteau Theorem imply that 
$$
\liminf_{n \rightarrow \infty} \pr \bigg( \sup_{ t \geq 0 } Y_n(t) > x \bigg) \geq 
\pr \big( \sup_{t \in [0,T] } {\mathcal Y}_{\infty}(t) > x \big).$$
Combining with the monotonicity of the supremum operator, and letting $T \rightarrow \infty$, it follows that for all $x \geq 0$,
\begin{equation}\label{secondparta}
\liminf_{n \rightarrow \infty} \pr \bigg( \sup_{ t \geq 0 } Y_n(t) > x \bigg) \geq \pr\big(\sup_{t \geq 0} {\mathcal Y}_{\infty}(t) > x \big).
\end{equation}
Combining (\ref{firstparta}) and (\ref{secondparta}), with the definition of weak convergence, completes the proof.
\endproof

\subsection{Proof of Theorem\ \ref{reed1}.}
In \cite{Reed.16}, Reed proves the following result.
\begin{theorem}[\cite{Reed.16}]\label{reed2}
Suppose that the HWR-$\alpha$ assumptions hold, and in addition S is deterministic (i.e. the system is $GI/D/n$).  Then $\lbrace n^{1 - \frac{1}{\alpha}} W^n_{A,S,B,\alpha}(\infty), n > B^{\frac{\alpha}{\alpha-1}} \rbrace$ converges in distribution to 
$\sup_{k \geq 0} \bigg( -(\frac{C_A}{C_{\alpha}})^{\frac{1}{\alpha}} \hat{S}_{\alpha,1}(k) - B k \bigg)$.
\end{theorem}
With Theorem\ \ref{reed2} in hand, we now apply the distributional Little's Law (and more generally the methodology of \cite{JMM.04}, which had previously been applied to the light-tailed setting) to derive the corresponding result for queue-lengths, Theorem\ \ref{reed1}.
\proof{Proof of Theorem\ \ref{reed1}}
Since the system is FCFS with i.i.d. inter-arrival and processing times, and processing times are determinstic (and hence there is no over-taking), the Distributional Little's Law applies (\cite{haji1971relation}), and we have 
\begin{equation}\label{dll1}
Q^n_{A,S,B,\alpha}(\infty) \sim A\bigg( \lambda_{n,B,\alpha} \big(1 + W^n_{A,S,B,\alpha}(\infty) \big) \bigg),
\end{equation}
with $A(t)_{t \geq 0}$ and $W^n_{A,S,B,\alpha}(\infty)$ independent.  Let $\lbrace A'_i, i \geq 1 \rbrace$ denote the sequence of inter-event times in ${\mathcal A}$ (i.e. corresponding to $A(t)_{t \geq 0}$), namely $A'_1$ is drawn from the equilibrium distribution, and $\lbrace A'_i, i \geq 2 \rbrace$ are i.i.d. distributed as $A$.  Then for all $x > 0$, $\pr \bigg( n^{-\frac{1}{\alpha}}\big( Q^n_{A,S,B,\alpha}(\infty) - n \big) \geq x \bigg)$ equals
\begin{eqnarray}
\ &\ &\ \pr\Bigg( A\bigg( \lambda_{n,B,\alpha} \big(1 + W^n_{A,S,B,\alpha}(\infty) \big) \bigg) \geq n + x n^{\frac{1}{\alpha}} \Bigg) \nonumber
\\&\ &\ \ \ =\ \ \ \pr\Bigg( \sum_{i=1}^{ \lceil n + x n^{\frac{1}{\alpha}} \rceil } A'_i \leq \lambda_{n,B,\alpha} \big(1 + W^n_{A,S,B,\alpha}(\infty) \big) \Bigg) \nonumber
\\&\ &\ \ \ =\ \ \  \pr\Bigg( \frac{\sum_{i=1}^{ \lceil n + x n^{\frac{1}{\alpha}} \rceil } (A'_i - 1)}{\lceil n + x n^{\frac{1}{\alpha}} \rceil^{\frac{1}{\alpha}}} \leq \frac{\lambda_{n,B,\alpha} \big(1 + W^n_{A,S,B,\alpha}(\infty) \big) - \lceil n + x n^{\frac{1}{\alpha}} \rceil}{ \lceil n + x n^{\frac{1}{\alpha}} \rceil^{\frac{1}{\alpha}} } \Bigg). \label{astay1}
\end{eqnarray}
It follows from Theorem\ \ref{genclt1} that 
\begin{equation}\label{astay2}
\bigg\lbrace \frac{\sum_{i=1}^{ \lceil n + x n^{\frac{1}{\alpha}} \rceil } (A'_i - 1)}{ \lceil n + x n^{\frac{1}{\alpha}} \rceil^{\frac{1}{\alpha}}} , n \geq 1 \bigg\rbrace\ \ \ \textrm{converges in distribution to}\ \ \ \big(\frac{C_A}{C_{\alpha}}\big)^{\frac{1}{\alpha}} S_{\alpha}(1,1,0).
\end{equation}
Theorem\ \ref{reed2} implies that 
\begin{equation}\label{astay3}
\bigg\lbrace \frac{\lambda_{n,B,\alpha} W^n_{A,S,B,\alpha}(\infty) }{ \lceil n + x n^{\frac{1}{\alpha}} \rceil^{\frac{1}{\alpha}} } , n \geq 1 \bigg\rbrace\ \ \ \textrm{converges in distribution to}\ \ \ 
\sup_{k \geq 0} \bigg( -(\frac{C_A}{C_{\alpha}})^{\frac{1}{\alpha}} \hat{S}_{\alpha,1}(k) - B k \bigg).
\end{equation}
Also, it is easily verified that
\begin{equation}\label{astay4}
\lim_{n \rightarrow \infty} \frac{ \lambda_{n,B,\alpha} - \lceil n + x n^\frac{1}{\alpha}\rceil }{ \lceil n + x n^{\frac{1}{\alpha}} \rceil^{\frac{1}{\alpha}} } = - B - x.
\end{equation}
As in \cite{JMM.04}, it then follows from the independence of $\lbrace A'_i, i \geq 1 \rbrace$ and $W^n_{A,S,B,\alpha}(\infty)$, and the CLT for triangular arrays (cf. \cite{Chung.74}) that for all $x$ which are continuity points of the c.d.f. of $\sup_{k \geq 1} \bigg( -(\frac{C_A}{C_{\alpha}})^{\frac{1}{\alpha}} \hat{S}_{\alpha,1}(k) - B k \bigg)$, it holds that
\begin{equation}\label{astay5}
\lim_{n \rightarrow \infty} \pr \bigg( n^{-\frac{1}{\alpha}}\big( Q^n_{A,S,B,\alpha}(\infty) - n \big) > x \bigg) = 
\pr\bigg( \sup_{k \geq 1} \bigg( -(\frac{C_A}{C_{\alpha}})^{\frac{1}{\alpha}} \hat{S}_{\alpha,1}(k) - B k \bigg) > x \bigg).
\end{equation}
The desired result then follows by applying the max-plus operator to both sides.
\endproof

\subsection{Proof of Corollary\ \ref{ldreed}.}
\proof{Proof of Corollary\ \ref{ldreed}}
Our approach is essentially identical to that used in \cite{MZ06}.  Let $X(t) \stackrel{\Delta}{=} -(\frac{C_A}{C_{\alpha}})^{\frac{1}{\alpha}} \hat{S}_{\alpha,1}(t) - B t$.  For $x > 0$, let $\tau(x) \stackrel{\Delta}{=} \inf\bigg\lbrace t \geq 0: X(t) \geq x \bigg\rbrace$, with $\tau(x) = \infty$ if the process never reaches a value greater than or equal to x.  In that case, for any $x > 0$ and $c \in (0,x)$, it follows from stationary and independent increments, and the strong Markov property, that 
\begin{eqnarray}
\pr\big( \sup_{t \geq 0} X(t) \geq x , \sup_{k \geq 0} X(k) \leq x - c \big) &\leq& \pr\bigg( \tau(x) < \infty, \inf_{s \in [\tau(x), \tau(x) + 1]} X(s) - X\big( \tau(x) \big) \leq - c \bigg) \nonumber
\\&=& \pr\big( \tau(x) < \infty\big) \times \pr\big( \inf_{s \in [0,1]} X(s) \leq - c \big)\nonumber
\\&=& \pr\big( \sup_{t \geq 0} X(t) \geq x \big) \pr\big( \inf_{s \in [0,1]} X(s) \leq - c \big).\label{uss1}
\end{eqnarray}
Combining with the fact that (by a union bound)
\begin{equation}\label{uss2}
\pr\big( \sup_{t \geq 0} X(t) \geq x \big) \leq \pr\big( \sup_{k \geq 0} X(k) > x - c\big) + \pr\big( \sup_{t \geq 0} X(t) \geq x , \sup_{k \geq 0} X(k) \leq x - c \big),
\end{equation}
we conclude that 
$$\pr\big( \sup_{t \geq 0} X(t) \geq x \big) \leq \pr\big( \sup_{k \geq 0} X(k) > x - c\big) + \pr\big( \sup_{t \geq 0} X(t) \geq x \big) \pr\big( \inf_{s \in [0,1]} X(s) \leq - c \big),$$
and thus
\begin{equation}\label{uss3}
\pr\big( \sup_{t \geq 0} X(t) \geq x \big) \leq \pr\big( \sup_{k \geq 0} X(k) > x - c\big) \times \bigg( \pr\big( \inf_{s \in [0,1]} X(s) > - c \big) \bigg)^{-1}.
\end{equation}
As $\inf_{s \in [0,1]} X(s)$ is a.s. finite, we may select c sufficiently large to ensure that $\pr\big( \inf_{s \in [0,1]} X(s) > - c \big) > 0$.  Then taking the appropriate limit as $x \rightarrow \infty$ (independent of the fixed value of $c$), and combining with Theorem\ \ref{beyondreedup1}, Corollary\ \ref{beyondreedupcor}, and Theorem\ \ref{reed1} completes the proof.
\endproof

\bibliographystyle{plainnat}

\end{document}